\documentclass[12pt]{amsart}
\usepackage{amsfonts,amssymb,color}
\usepackage[mathscr]{eucal}
\usepackage{amsmath, amsthm}
\usepackage{mathrsfs}
\usepackage{amsbsy}
\usepackage{dsfont}
\usepackage{bbm}
\usepackage{wasysym}
\usepackage{stmaryrd}
\input xypic
\xyoption{all}
\begin{document}
\baselineskip = 5,2mm
\newcommand \ZZ {{\mathbb Z}} 
\newcommand \FF {{\mathbb F}} %
\newcommand \NN {{\mathbb N}} 
\newcommand \QQ {{\mathbb Q}} 
\newcommand \RR {{\mathbb R}} 
\newcommand \CC {{\mathbb C}} 
\newcommand \PR {{\mathbb P}} 
\newcommand \AF {{\mathbb A}} 
\newcommand \bcA {{\mathscr A}}
\newcommand \bcB {{\mathscr B}}
\newcommand \bcC {{\mathscr C}}
\newcommand \bcD {{\mathscr D}}
\newcommand \bcE {{\mathscr E}}
\newcommand \bcF {{\mathscr F}}
\newcommand \bcG {{\mathscr G}}
\newcommand \bcH {{\mathscr H}}
\newcommand \bcM {{\mathscr M}}
\newcommand \bcN {{\mathscr N}}
\newcommand \bcI {{\mathscr I}}
\newcommand \bcJ {{\mathscr J}}
\newcommand \bcK {{\mathscr K}}
\newcommand \bcL {{\mathscr L}}
\newcommand \bcO {{\mathscr O}}
\newcommand \bcP {{\mathscr P}}
\newcommand \bcQ {{\mathscr Q}}
\newcommand \bcR {{\mathscr R}}
\newcommand \bcS {{\mathscr S}}
\newcommand \bcT {{\mathscr T}}
\newcommand \bcU {{\mathscr U}}
\newcommand \bcV {{\mathscr V}}
\newcommand \bcW {{\mathscr W}}
\newcommand \bcX {{\mathscr X}}
\newcommand \bcY {{\mathscr Y}}
\newcommand \bcZ {{\mathscr Z}}
\newcommand \goa {{\mathfrak a}}
\newcommand \gob {{\mathfrak b}}
\newcommand \goc {{\mathfrak c}}
\newcommand \gom {{\mathfrak m}}
\newcommand \gop {{\mathfrak p}}
\newcommand \goT {{\mathfrak T}}
\newcommand \goC {{\mathfrak C}}
\newcommand \goD {{\mathfrak D}}
\newcommand \goM {{\mathfrak M}}
\newcommand \goN {{\mathfrak N}}
\newcommand \goS {{\mathfrak S}}
\newcommand \goH {{\mathfrak H}}
\newcommand \Noe {{\sf N}}
\newcommand \Noesn {{{\sf N}_0}}
\newcommand \uno {{\mathbbm 1}}
\newcommand \Le {{\mathbbm L}}
\newcommand \Spec {{\rm {Spec}}}
\newcommand \Pic {{\rm {Pic}}}
\newcommand \Jac {{\rm {Jac}}}
\newcommand \Alb {{\rm {Alb}}}
\newcommand \NS {{{NS}}}
\newcommand \Corr {{Corr}}
\newcommand \Chow {{\mathscr C}}
\newcommand \Sym {{\rm {Sym}}}
\newcommand \Alt {{\rm {Alt}}}
\newcommand \Prym {{\rm {Prym}}}
\newcommand \cone {{\rm {cone}}}
\newcommand \cha {{\rm {char}}}
\newcommand \eff {{\rm {eff}}}
\newcommand \rat {{\rm rat}}
\newcommand \tr {{\rm {tr}}}
\newcommand \pr {{\rm {pr}}}
\newcommand \ev {{\it {ev}}}
\newcommand \interior {{\rm {Int}}}
\newcommand \sep {{\rm {sep}}}
\newcommand \td {{\rm {tdeg}}}
\newcommand \alg {{\rm {alg}}}
\newcommand \im {{\rm im}}
\newcommand \et {\rm {\acute e t}}
\newcommand \Zar {\rm {Zar}}
\newcommand \Nis {\rm {Nis}}
\newcommand \cdh {\rm {cdh}}
\newcommand \h {\rm {h}}
\newcommand \Pre {{\mathscr P}}
\newcommand \Shv {{\sf Shv}}
\newcommand \Funct {{\rm Funct}}
\newcommand \op {{\rm op}}
\newcommand \Hom {{\rm Hom}}
\newcommand \uHom {{\underline {\rm Hom}}}
\newcommand \Hilb {{\rm Hilb}}
\newcommand \Sch {{\mathscr S\! }{\it ch}}
\newcommand \cHilb {{\mathscr H\! }{\it ilb}}
\newcommand \cHom {{\mathscr H\! }{\it om}}
\newcommand \cExt {{\mathscr E\! }{\it xt}}
\newcommand \colim {{{\rm colim}\, }} 
\newcommand \End {{\rm {End}}}
\newcommand \coker {{\rm {coker}}}
\newcommand \id {{\rm {id}}}
\newcommand \van {{\rm {van}}}
\newcommand \spc {{\rm {sp}}}
\newcommand \Ob {{\rm Ob}}
\newcommand \Aut {{\rm Aut}}
\newcommand \cor {{\rm {cor}}}
\newcommand \res {{\rm {res}}}
\newcommand \Gal {{\rm {Gal}}}
\newcommand \PGL {{\rm {PGL}}}
\newcommand \Gr {{\rm {Gr}}}
\newcommand \Bl {{\rm {Bl}}}
\newcommand \Sing {{\rm {Sing}}}
\newcommand \spn {{\rm {span}}}
\newcommand \Nm {{\rm {Nm}}}
\newcommand \inv {{\rm {inv}}}
\newcommand \codim {{\rm {codim}}}
\newcommand \ptr {{\pi _2^{\rm tr}}}
\newcommand \sg {{\Sigma }}
\newcommand \CHM {{\mathscr C\! \mathscr M}}
\newcommand \DM {{\sf DM}}
\newcommand \FS {{FS}}
\newcommand \MM {{\mathscr M\! \mathscr M}}
\newcommand \HS {{\mathscr H\! \mathscr S}}
\newcommand \MHS {{\mathscr M\! \mathscr H\! \mathscr S}}
\newcommand \Vect {{\mathscr V\! ect}} %
\newcommand \Gm {{{\mathbb G}_{\rm m}}}
\newcommand \trdeg {{\rm {tr.deg}}}
\newcommand \Mon {{\rm Mon }}
\newcommand \jacob {\tiny {\wasylozenge }}
\newcommand \univ {\tiny {\wasylozenge }}
\newcommand \tame {\rm {tame }}
\newcommand \prym {\tiny {\Bowtie }}
\newcommand \znak {{\natural }}
\newcommand \qand {{\quad \hbox{and}\quad }}
\newcommand \qqand {{\quad \quad \hbox{and}\quad \quad }}
\newcommand \heither {{\hbox{either}\quad }}
\newcommand \qor {{\quad \hbox{or}\quad }}
\newtheorem{theorem}{Theorem}
\newtheorem{lemma}[theorem]{Lemma}
\newtheorem{corollary}[theorem]{Corollary}
\newtheorem{proposition}[theorem]{Proposition}
\newtheorem{remark}[theorem]{Remark}
\newtheorem{definition}[theorem]{Definition}
\newtheorem{conjecture}[theorem]{Conjecture}
\newtheorem{example}[theorem]{Example}
\newtheorem{question}[theorem]{Question}
\newtheorem{assumption}[theorem]{Assumption}
\newtheorem{fact}[theorem]{Fact}
\newtheorem{crucialquestion}[theorem]{Crucial Question}
\newcommand \lra {\longrightarrow}
\newcommand \hra {\hookrightarrow}
\def\blue {\color{blue}}
\def\red {\color{red}}
\def\green {\color{green}}
\newenvironment{pf}{\par\noindent{\em Proof}.}{\hfill\framebox(6,6)
\par\medskip}
\title[Rational equivalence for 1-cycles on cubic hypersurfaces in $\PR ^5$]
{\bf \'Etale monodromy and rational equivalence for $1$-cycles on cubic hypersurfaces in $\PR ^5$}
\author{Kalyan Banerjee and Vladimir Guletski\u \i }

\date{28 March 2018}


\begin{abstract}
\noindent Let $k$ be an uncountable algebraically closed field of characteristic $0$, and let $X$ be a smooth projective connected variety of dimension $2p$, appropriately embedded into $\PR ^m$ over $k$. Let $Y$ be a hyperplane section of $X$, and let $A^p(Y)$ and $A^{p+1}(X)$ be the groups of algebraically trivial algebraic cycles of codimension $p$ and $p+1$ modulo rational equivalence on $Y$ and $X$ respectively. Assume that, whenever $Y$ is smooth, the group $A^p(Y)$ is regularly parametrized by an abelian variety $A$ and coincides with the subgroup of degree $0$ classes in the Chow group $CH^p(Y)$. In the paper we prove that the kernel of the push-forward homomorphism from $A^p(Y)$ to $A^{p+1}(X)$ is the union of a countable collection of shifts of a certain abelian subvariety $A_0$ inside $A$. For a very general section $Y$ either $A_0=0$ or $A_0$ coincides with an abelian subvariety $A_1$ in $A$ whose tangent space is the group of vanishing cycles $H^{2p-1}(Y)_{\van }$. Then we apply these general results to sections of a smooth cubic fourfold in $\PR ^5$.
\end{abstract}
















\maketitle

\tableofcontents

\section{Introduction}
\label{intro}

Let $X$ be a smooth projective variety over an algebraically closed field. The Picard-Lefschetz theory yields that the monodromy action on the $(n-1)$-th vanishing cohomology of a smooth section of the variety $X$ is irreducible. The proof of this fact in terms of \'etale cohomology groups was given by N. Katz in his second article (Expos\'e XVIII) in \cite{SGA7II}, and by P. Deligne in \cite{WeilConjI}. It is also explained in terms of singular cohomology in Section 7.3 of Lamotke's paper \cite{Lamotke}.

The irreducibility of the monodromy action plays an important role in the Hodge theoretical study of algebraic varieties over $\CC $, and it was amply utilized in the work by C. Voisin, see, for example, \cite{Voisin: Theoreme de Torelli}, \cite{Voisin: Sur les 0-cycles}, \cite{Voisin: Variations}, \cite{Voisin: The generalized Hodge and Bloch I} and \cite{Voisin: Universal CH_0}. More importantly, the irreducibility of monodromy action on cohomology affects algebraic cycles through Hodge theory, see pages 304 - 305 in the second volume of the book \cite{Voisin: The Book}, and Proposition 2.4 on page 854 in \cite{Voisin: Sympl. involutions}.

To explain the latter idea, let $X$ be a smooth projective complex surface, embedded into a projective space, let $Y$ be a general hyperplane section of $X$ with the Jacobian $A=\Jac (Y)$, and let $A_0(X)$ be the Chow group of $0$-cycles of degree $0$ on $X$. Let, furthermore, $A_1$ be an abelian subvariety in $A$ which corresponds, via Hodge theory, to the vanishing cohomology in $H^1(Y)$. Then the kernel of the push-forward homomorphism from $A$ to $A_0(X)$ is a countable union of shifts of a certain abelian subvariety $A_0$ inside $A_1$. The cohomological monodromy argument implies that, for a general $Y$, either $A_0=A_1$ or $A_0=0$, see the top of page 305 in the second volume of \cite{Voisin: The Book}. Clearly, this alternative for $A_0$ must play an important role in the study of $0$-cycles on surfaces, especially in the context of Bloch's conjecture.


The main objective of the present paper is to initiate a systematic study of Voisin's cycle-theoretic monodromy argument in as much generality as possible. To be more precise, we will apply the concept to algebraic cycles of dimension $p-1$ on a smooth projective variety $X$ of even dimension $2p$, embedded into a projective space, such that, if $Y$ is a smooth hyperplane section of $X$, the ``continuous" Chow group $A^p(Y)$ is weakly representable in the sense of Spencer Bloch, see Definition 1.1 in \cite{BlochAnExample} or Definition 3.3. in \cite{BlochMurre}.

Our second major objective is to develop a cycle-theoretic monodromy argument working over an abstract field of definition. Though we were not able to avoid the uncountability and $0$ characteristic of the ground field in the present paper, we strongly believe that these two requirements can be effectively omitted in a more subtle arithmetic theory, in which the intrinsic nature of the abelian variety $A_0$ will be revealed.

Let us now describe the results of the paper. Let $k$ be an uncountable algebraically closed field of characteristic $0$, and let
  $$
  r:Y\to X
  $$
be a codimension $e$ closed embedding of smooth projective varieties over $k$. The morphism $r$ induces a push-forward homomorphism
  $$
  r_*:A^p(Y)\to A^{p+e}(X)
  $$
on the Chow groups of algebraically trivial algebraic cycles, and our aim is to study the kernel
  $$
  K=\ker (r_*)
  $$
under the assumption that the group $A^p(Y)$ is weakly representable, or, equivalently, it can be regularly parametrized by an abelian variety $A$ over $k$.

Fix an embedding of the variety $X$ into a projective space $\PR ^m$, such that $X$ is not contained in a hyperplane. For simplicity, and as it is sufficient for applications we have in mind, we will also assume that the group $CH^p(Y)_{\deg =0}$, defined by the embedding of $Y$ into $\PR ^m$, coincides with the group $A^p(Y)$. Clearly, this assumption is satisfied if $p=1$ and the N\'eron-Severi group is of rank $1$, or if $p$ is the dimension of $Y$ (i.e. we study $0$-cycles). They are also satisfied when $Y$ is a Fano threefold of Picard number $1$ inside $X$ and $p=2$. For example, when $Y$ is a smooth section of a smooth cubic fourfold in $\PR ^5$.

Our first result is a generalization of the Mumford-Roitman countability lemma for $0$-cycles, see \cite{Mumford} and \cite{Roitman}, to cycles of positive dimension.

\bigskip

\begin{itemize}

\item[]{}
{\rm THEOREM A.} {\it There exists an abelian subvariety $A_0$ in $A$, and a countable set $\Xi $ of closed points on $A$, such that
  $$
  K=\cup _{x\in \Xi }(x+A_0)
  $$
inside the abelian variety $A$.
}
\end{itemize}

\bigskip

The presentation of the group $A^p(Y)$ by the abelian variety $A$ provides a homomorphism from $H^1(A)$ to $H^{2p-1}(Y)$, in terms of $l$-adic cohomology groups. Assuming that this homomorphism is an isomorphism, which is known to be true for $p\leq 2$, and also using the Tate conjecture for abelian varieties proven by Faltings, one can construct an abelian subvariety $A_1$ in $A$ whose tangent space is the kernel of the Gysin homomorphism from $H^{2p-1}(Y)$ to $H^{2(p+e)-1}(X)$. We prove in the paper that the abelian variety $A_0$ is a subvariety in $A_1$.

Let now $S$ be an integral algebraic scheme over $k$, let $\bar \eta $ be the geometric generic point of $S$, and choose a $c$-open subset $U$ in $S$, such that for any closed point $P$ in $U$ there is a scheme-theoretical isomorphism between $\bar \eta $ and $P$ over $\QQ $. Consider a closed embedding $\bcY \to \bcX $ of smooth families over $S$. The scheme-theoretic isomorphisms $\bar \eta \simeq P$ induce scheme-theoretic isomorphisms $\varkappa _P$ between $\bcY _P$ and $\bcY _{\bar \eta }$ over $\QQ $. Assume that $A^p(\bcY _{\bar \eta })$ is presented by an abelian variety $A_{\bar \eta }$ and $A^p(\bcY _P)$ is presented by an abelian variety $A_P$, for every closed point $P$ in $U$. Then the isomorphisms $\varkappa _P$ induce isomorphisms $\kappa _P$ between $A_P$ and $A_{\bar \eta }$ compatible with the isomorphisms on Chow groups induced by the isomorphisms $\varkappa _P$. In the paper we show that $\kappa _P(A_{P,0})=A_{\bar \eta ,0}$ and $\kappa _P(A_{P,1})=A_{\bar \eta ,1}$ for every $k$-point $P$ in $U$. In other words, one can study the varieties $A_0$ in a family either working at the geometric generic point or at a very general closed point on the base scheme $S$. The $c$-open set $U$ is not unique, of course, and all remains the same over any of them. What happens to $A_0$ beyond the union of such $c$-open sets, is a big and important question which deserves a separate research programme.

Within this paper we are interested in the case where the family in question is a family of smooth hyperplane sections, so that we can enhance the study of the abelian variety $A_0$ by the monodromy action. So let $X$ be a smooth projective variety of even dimension $2p$ over the ground field $k$, embedded into a projective space $\PR ^m$, let $\bcX =X\times {\PR ^m}^{\vee }$, where ${\PR ^m}^{\vee }$ is the dual projective space, and let $\bcY $ be the intersection of $\bcX $ with the universal hyperplane inside $\PR ^m\times {\PR ^m}^{\vee }$. Let also $T$ be the complement to the discriminant locus inside the dual projective space, and consider the family $\bcY _T\to T$ of smooth hyperplane sections of the variety $X$ over $k$. Clearly, $\bcY _T$ is embedded into $\bcX _T$ over $T$. Let, furthermore, $\xi $ be the generic point of $T$, let $\bar \xi $ be the corresponding geometric generic point, and choose a $c$-open subset $U$ in $T$, so that $k$-points $P$ in $U$ are scheme-theoretically isomorphic to $\bar \xi $. Now again, assuming our standard assumptions for the fibres $Y_{\bar \xi }$ and $Y_P$, and applying the $l$-adic monodromy argument in Lefschetz pencils passing through $U$, in conjunction with Theorem A, we obtain

\bigskip

\begin{itemize}

\item[]{}
{\rm THEOREM B.} {\it In terms above, either $A_{\bar \xi ,0}=0$ or $A_{\bar \xi ,0}=A_{\bar \xi ,1}$. Respectively, either $A_{P,0}=0$ or $A_{P,0}=A_{P,1}$, for any closed point $P$ in $U$.
}
\end{itemize}

\bigskip

Notice that the assumptions of Theorem B are satisfied, for example, for all smooth projective fourfolds $X$ whose very general hyperplane sections are Fano varieties of Picard number $1$. In all such cases Theorem B brings new information about rational equivalence of algebraic $1$-cycles on the fourfold $X$.

Let now $\tilde Y_P$ be the resolution of singularities on $Y_P$, and set $\tilde Y_P=Y_P$ whenever the section $Y_P$ is smooth. Assume that $p\leq 2$ and the group $A^p(\tilde Y_P)$ is weakly representable for every section $Y_P$ having at worst one ordinary double point. The next theorem generalizes Voisin's result on page 305 in Volume II of \cite{Voisin: The Book}.

\bigskip

\begin{itemize}

\item[]{}
{\rm THEOREM C.} {\it If the group $A^{p+1}(X)$ is not rationally weakly representable, it follows that the kernel of the push-forward homomorphism from $A^p(Y_P)$ to $A^{p+1}(X)$ is countable, for a very general hyperplane section $Y_P$.}

\end{itemize}

\bigskip

The practical meaning of Theorem C is as follows. Suppose we want to prove that $A^{p+1}(X)$ is weakly representable. Then, by Theorem C, ``all we need" is to find a positive-dimensional variety in the kernel of the homomorphism from $A^p(Y)$ to $A^{p+1}(X)$, for a very general ample section $Y$ on $X$.

Our original motivation for proving Theorem C was to understand more about the structure of the huge Chow group $CH^3(X)$ for a general cubic hypersurface $X$ in $\PR ^5$, whence the title of the paper. Recall that $CH^3(X)$ is generated by lines, see Theorem 1.1 in \cite{Shen_relations}, and $A^3(X)$ is not weakly representable by Theorem 0.5 in \cite{Schoen}. Theorem C tells us that one can think of $A^3(X)$ as a collection of Prymians of smooth hyperplane sections modulo countable kernels generated by $1$-cycles rationally equivalent to $0$ on $X$. We expect that these countable kernels are of deep arithmetical nature relevant to the famous non-rationality conjecture.

\medskip

{\sc Acknowledgements.} We are grateful to Sergey Gorchinskiy, Mingmin Shen, Alexander Tikhomirov, Claire Voisin and Yuri Zarhin for useful discussions relevant to the theme of this paper. The second named author has been partially supported by EPSRC grant EP/I034017/1. The first named author acknowledges EPSRC GTA PhD studentship in Liverpool.

\section{Problem setting and standard assumptions}
\label{setting&assumptions}

The purpose of the first section is to set up the main problem and justify three basic assumptions which we will keep throughout the paper.

Let $k$ be an algebraically closed field of characteristic $0$. For an algebraic variety $Y$ over $k$, let $A^p(Y)$ be the subgroup in $CH^p(Y)$ generated by algebraically trivial algebraic cycles on $Y$. Suppose $V$ is another algebraic variety over $F$ with a fixed closed point $P_0$ on it. Let $Z$ be an algebraic cycle of codimension $n$ on the product $V\times Y$. For any closed point $P$ on $V$ we have the flat pullback $Z(P)$ of the cycle $Z$ to the fibre $Y_P=Y$ of the projection $V\times Y\to V$ at $P$. Then $Z(P)$ is a codimension $n$ algebraic cycles on the variety $Y$. The cycle
  $$
  Z(P)-Z(P_0)
  $$
is algebraically trivial on $Y$, and we obtain a map
  $$
  V\to A^p(Y)
  $$
sending any closed point $P$ on $V$ to the class of the cycle $Z(P)-Z(P_0)$ on $Y$. This map is nothing but the algebraic family of codimension $n$ algebraically trivial cycle classes on $Y$ determined by the algebraic cycle $Z$ on $V\times Y$ and the point $P_0$ on the parameter variety $V$. The next definition appears in Murre's paper \cite{Murre2}, and is important for what follows.

\begin{definition}
\label{regularmaps}
{\rm Let $A$ be an abelian variety over $k$. A group homomorphism $A^p(Y)\to A$ is said to be {\it regular} if its pre-composition with any family of algebraic cycles $V\to A^p(Y)$ in the sense above is a regular morphism over $k$. A regular homomorphism
	$$
	\psi :A^p(Y)\to A
	$$
to an abelian variety $A$ over $k$ is said to be {\it universal} if, having another regular homomorphism $A^p(Y)\to B$, there exists a unique homomorphism of abelian varieties $A\to B$, such that the obvious diagram commutes, see \cite{Murre2}, page 981. Clearly, if $\psi $ exists, then it is an epimorphism in the category of abelian groups.
}
\end{definition}



Let
  $$
  r:Y\to X
  $$
be a closed embedding of smooth projective connected varieties over $k$ of codimension $e$, let
  $$
  r_*:A^p(Y)\to A^{p+e}(X)
  $$
be the push-forward homomorphism induced by the proper morphism $r$, and let $K$ be the kernel of the homomorphism $r_*$. Our aim is to study the kernel $K$ of the push-forward homomorphism $r_*$. Certainly, it is difficult to study $K$ in general, and therefore we need to impose some reasonable assumptions on the Chow group $A^p(Y)$.

\bigskip

\noindent {\it Assumption 1}

\medskip

\noindent Our first assumption is that the universal regular epimorphism
  $$
  \psi :A^p(Y)\to A
  $$
exists, and that the group $A^p(Y)$ is weakly representable in the sense of Spencer Bloch, in which case the universal regular homomorphism $\psi $ is an isomorphism, so that we can identify $A^p(Y)$ and $A$ by means of $\psi $, see \cite{BlochAnExample} or \cite{BlochMurre}.

\begin{remark}
{\rm Clearly, the universal homomorphism $\psi $ exists and is an isomorphism if $p=1$. The main result of \cite{Murre2} asserts that $\psi $ exists in case $p=2$. Therefore, Assumption 1 is satisfied whenever $p=2$ and $A^2(Y)$ is weakly representable in the sense of Bloch, see \cite{BlochAnExample} and \cite{BlochMurre}. Notice that if the group $A^3(Y)$ of $0$-cycles on $Y$ is weakly representable, so is the group $A^2(Y)$, see Lemma 3.1 in \cite{GorchGulFano}. Therefore, Assumption 1 is satisfied whenever $A^3(Y)$ is representable. In particular, it is satisfied when $Y$ is rationally connected. This is so, for example, if $Y$ is a Fano threefold inside a smooth projective variety $X$ over $k$.
}
\end{remark}

\bigskip

\noindent {\it Assumption 2}

\medskip

\noindent Fix an embedding of the variety $X$ into a projective space $\PR ^m$. Since $Y$ is a subvariety in $X$, it induces the embedding of $Y$ into the same space $\PR ^m$. Obviously, $A^p(Y)\subset CH^p(Y)_{\deg =0}$. For simplicity, and as it is sufficient for applications we have in mind, we will also assume that
  $$
  A^p(Y)=CH^p(Y)_{\deg =0}\; .
  $$

\begin{remark}
{\rm If $p=1$, this assumption is obviously satisfied. If $p=2$, Assumption 2 is satisfied, for example, for all Fano threefolds $Y$ inside a smooth projective fourfold $X$ over $k$, such that the Picard number of $Y$ is equal to $1$. In particular, the second assumption is satisfied for smooth sections $Y$ of a smooth cubic fourfold in $\PR ^4$.
}
\end{remark}

\bigskip

\noindent {\it Assumption 3}

\medskip

\noindent Since $\psi $ is an isomorphism, and hence the group $A^p(Y)$ is weakly representable, there exists a smooth projective curve $\Gamma$, a cycle $Z$ of codimension $p$ on $\Gamma \times Y$, and an algebraic subgroup $G\subset J_{\Gamma }$ in the Jacobian variety $J_{\Gamma }$, such that the induced homomorphism
  $$
  z_*:J_{\Gamma }=A^1(\Gamma )\to A^p(Y)\simeq A
  $$
is surjective, and its kernel is the group $G$. Here $z$ is the cycle class of $Z$ in the Chow group $CH^p(\Gamma \times Y)$. Furthermore, the class $z$ gives the morphism
  $$
  z:M(\Gamma )\otimes \Le ^{p-1}\to M(Y)\; ,
  $$
where $M(-)$ is the functor from smooth projective varieties over $k$ to (contravariant) Chow motives over $k$, $\Le $ is the Lefschetz motive and $\Le ^n$ is the $n$-fold tensor power of $\Le $. Fix a point on the curve $\Gamma $ and consider the induced embedding
  $$
  i_{\Gamma }:\Gamma \to J_{\Gamma }\; .
  $$
Let also
  $$
  \alpha :J_{\Gamma }\to A
  $$
be the projection from the Jacobian $J_{\Gamma }$ onto the abelian variety $A=J_{\Gamma }/G$. Define $w$ to be the composition
  $$
  z\circ (M(\alpha \circ i_{\Gamma })\otimes \id _{\Le ^{p-1}})
  $$
in the category of Chow motives with coefficients in $\ZZ $. Then $w$ is a morphism from the motive $M(A)\otimes \Le ^{p-1}$ to the motive $M(Y)$, which induces the homomorphism
  \begin{equation}
  \label{scebyaraki}
  w_*:H^1_{\et }(A,\QQ _l(1-p))\to
  H^{2p-1}_{\et }(Y,\QQ _l)\; ,
  \end{equation}
see \cite{GorchGulFano}. Our third assumption is that the homomorphism $w_*$ is an isomorphism of cohomology groups.

\begin{remark}
\label{lubanskiiles}
{\rm If $p=1$ and $\dim (Y)=1$, then $w_*$ is an isomorphism by the standard argument. If $p=2$ and $\dim (Y)=3$, Assumption 3 is satisfied by Lemma 4.3 in \cite{GorchGulFano}. Conjecturally, $w_*$ is an isomorphism for any $p>2$ as well, but we did not prove that in the paper. The reason for that is that Lemma 4.3 from \cite{GorchGulFano} uses some result by Merkurjev and Suslin on the injectivity of Bloch's map $\lambda^2_l$, see \cite{MS}. Thought we believe that the Bloch-Kato conjecture, which is now a theorem due to Voevodsky and Rost, can help us to prove that $w_*$ is an isomorphism for an arbitrary $p$, this may well be quite a big piece of work, deserving a separate paper to be worked it out.

Of course, if $k=\CC $, then the isomorphism between $H^1(A,\CC )$ and $H^{2p-1}(Y,\CC )$ can be easily achieved by the Hodge-theoretical methods, and Assumption 3 is always satisfied.
}
\end{remark}

Consider the Gysin homomorphism
  \begin{equation}
  \label{syadziba}
  H^{2p-1}_{\et }(Y,\QQ _l)
  \stackrel{r_*}{\lra }
  H^{2(p+e)-1}_{\et }(X,\QQ _l)\; .
  \end{equation}
induced by the closed embedding $r$ on the $l$-adic cohomology groups. Assumption 3 gives us an advantage that we now can describe the kernel of the homomorphism (\ref{syadziba}) in terms of an abelian subvariety in $A$, transporting Hodge-theoretical arguments in to $l$-adic representations.

Indeed, let
  $$
  A_{l^n}=\ker (A\stackrel{l^n}{\to }A)
  $$
be the $l^n$-torsion subgroup in the abelian variety $A$ over $k_0$, let
  $$
  T_l(A)=\lim A_{l^n}
  $$
be the Tate module of $A$, and let
  $$
  V_l(A)=T_l(A)\otimes _{\ZZ _l}\QQ _l\; .
  $$

Let $G$ be the image in the group $H^1_{\et }(A,\QQ _l(1-p))$ of the kernel of the Gysin homomorphism (\ref{syadziba}) under the isomorphism $w_*^{-1}$, inverse to (\ref{scebyaraki}). Then $G$ induces a $\QQ _l$-vector subspace in $H^1_{\et }(A,\QQ _l)$. But the group $H^1_{\et }(A,\QQ _l)$ is isomorphic to the dual vector space $\Hom _{\QQ _l}(V_l(A),\QQ _l)$. Therefore, the kernel $G$ induces a vector subaspace $G'$ in $\Hom _{\QQ _l}(V_l(A),\QQ _l)$. Since the space $V_l(A)$ is finite-dimensional, the dual space $\Hom _{\QQ _l}(V_l(A),\QQ _l)$ is isomorphic to $V_l(A)$, and therefore $G'$ induces a vector subspace in $V_l(A)$. The latter vector subspace determines an idempotent $e_G$ in the associative ring $\End _{\QQ _l}(V_l(A))$.

Now, without loss of generality, one can temporarily assume that $k$ is the algebraic closure of a field which is finitely generated over $\QQ $. In such a case, the Tate conjecture for abelian varieties, proven by Faltings, tells us that the canonical $l$-adic representation
  $$
  \End (A)\otimes \QQ _l\to \End (V_l(A))
  $$
is an isomorphism, see the article \cite{Faltings} (or page 72 in \cite{Tate}, or page 74 in \cite{Shatz}). Therefore, the idempotent $e_G$ induces an idempotent $\tilde e_G$ in the associative ring $\End (A)\otimes \QQ _l$. This idempotent determines a unique, up to an isogeny, abelian subvariety
  $$
  A_1\subset A
  $$
in the abelian variety $A$, such that the image of the injective homomorphism
  $$
  H^1_{\acute e t}(A_1,\QQ _l)(1-p)\to
  H^1_{\acute e t}(A,\QQ _l)(1-p)\; ,
  $$
induced by the inclusion $A_1\subset A$, coincides with the kernel of the composition of the isomorphism $w_*$ with the Gysin homomorphism $r_*$ from $H^{2p-1}_{\acute e t}(Y,\QQ _l)$ to $H^{2p+1}_{\acute e t}(X,\QQ _l)$.

\begin{remark}
{\rm If $p=1$, $\dim (Y)=1$ and $\dim (X)=2$, then $A_1$ is the connected component of the kernel of the induced homomorphism from the Jacobian $A$ of the curve $Y$ to the Albenese variety $\Alb (X)$ of the surface $X$.
}
\end{remark}

\begin{remark}
{\rm In the applications below we will be dealing with the case when $X$ is embedded into a projective space, the dimension of $X$ is $2p$ and $Y$ is a smooth hyperplane section of $X$, so that $\dim (Y)=2p-1$ and $e=1$. Assume also that $p=1$ or $2$, in order to have that $w_*$ is an isomorphism by Remark \ref{lubanskiiles}. If $H^{2p+1}_{\et }(X,\QQ _l)$ vanishes, then the primitive cohomology group coincides with $H^{2p-1}_{\et }(Y,\QQ _l)$, in which case $A_1=A$. This is so, for example, when $X$ is a smooth hypersurface in $\PR ^{2p+1}$. If the group does not $H^{2p+1}_{\et }(X,\QQ _l)$ vanish, the abelian variety $A_1$ can be smaller than the variety $A$.
}
\end{remark}

\begin{remark}
{\rm If, moreover, $k=\CC $, then $A_1$ can be described Hodge-theoretically. Indeed, for any algebraic variety $V$ over $\CC $ and any non-negative integer $n$ the \'etale cohomology group $H^n_{\et }(V,\QQ _l)$ is functorially isomorphic to the singular cohomology group $H^*(V(\CC ),\QQ _l)=H^*(V(\CC ),\QQ )\otimes \QQ _l$. The \'etale cohomology groups with coefficients in $\QQ _l$ can be further tensored with the algebraic closure $\bar \QQ _l$ of the $l$-adic field over $\QQ _l$. Fixing an isomorphism between $\bar \QQ _l$ and $\CC $, the \'etale cohomology groups $H^*_{\et }(-,\bar \QQ _l)$ are functorially isomorphic to the singular cohomology groups $H^*(-,\CC )$. The Gysin homomorphism $r_*$ from $H^{2p-1}(Y,\CC )$ to $H^{2p+1}(X,\CC )$ is a morphism of Hodge structures, so that its kernel $H_1$ is a Hodge substructure in $H^{2p-1}(Y,\CC )$. Suppose $p=2$. By Remark \ref{lubanskiiles} the group $H^3(Y,\CC )$ is isomorphic to the group $H^1(A,\CC )$ via the homomorphism $w_*$, and $w_*$ is obviously a morphism of Hodge structures too. It follows that $H_1$ is of weight $1$. This gives the abelian subvariety $A_1$ in $A$, where $A=J^2(Y)_{\alg }$ is the intermedian Jacobian of the threefold $Y$ (see \cite{Murre2}).
}
\end{remark}

\section{The generalization of the Mumford-Roitman theorem}
\label{Chowmonoids}

In this section we generalize a certain result due to Mumford and Roitman, appeared first in \cite{Mumford} and then in \cite{Roitman}, to algebraic cycles of positive dimension. For that purpose we shall use the theory of relative cycles developed by Suslin and Voevodsky in \cite{SuslinVoevodskyChow} and, independently, by Koll\'ar in \cite{Kollar}.

So, let $\Noe $ be the category of locally Noetherian schemes over $k$. Let $X$ be a scheme of finite type over $k$ and consider the presheaf of monoids\footnote{all monoids in this paper will be commutative monoids} $\bcC _r(X)$ on $\Noe $, where for any scheme $S$ in $\Noe $ the value $\bcC _r(X)(S)$ is the monoid $\bcC _r(X\times S/S)$ freely generated by relative cycles on $X\times S$ of relative dimension $r$ over $S$, in the sense of Suslin and Voevodsky, and the pullbacks are the pullbacks as constructed in Section 3 in \cite{SuslinVoevodskyChow}. To understand why the monoids $\bcC _r(X)(S)$ are free, see Corollary 3.1.6, or Corollaries 3.4.5 and 3.4.6, in loc.cit. The presheaf $\bcC _r(X)$ is actually a sheaf in $\cdh $, and hence in the Nisnevich topology on $\Noe $, see Theorem 4.2.9 in \cite{SuslinVoevodskyChow}. If $X$ is equidimensional, we will also write $\bcC ^p(X)(S)$, or $\bcC ^p(X\times S/S)$, for the same monoids of relative cycles of relative codimension $p$, where $p=\dim (X)-r$. If $X$ is projective over $k$, we fix a closed embedding of $X$ into a projective space $\PR ^m$ over $k$ and consider the subsheaf $\bcC ^p_d(X)$ of relative cycles of degree $d$ in $\bcC ^p(X)$, whose sections $\bcC ^p_d(X\times S/S)$ are generated by relative cycles of relative codimension $p$ and degree $d$, where the degree is understood with regard to the closed embedding of $X$ into $\PR ^m$.

Since now we will always assume that $X$ is equidimensional and projective over $k$, and that the embedding of $X$ into $\PR ^m$ is fixed. Corollary 4.4.13 in \cite{SuslinVoevodskyChow} says then that the sheafification of the presheaf $\bcC ^p_d(X)$ in $h$-topology on the category $\Noe $ is representable by a {\it Chow scheme} $C^p_d(X)$, projective over $k$. Let
  $$
  C^p(X)=\coprod _{d\geq 0}C^p_d(X)
  $$
be the {\it total Chow scheme}, where the coproduct is taken in the category $\Noe $.

Now, let $\Noesn $ be the full subcategory of seminormal schemes in $\Noe $. Since $k$ is of characteristic $0$, the above $h$-representability can be replaced by the usual representability, if we restrict the presheaves on the category $\Noesn $, see Theorem 3.21 in \cite{Kollar}, and notice that the Suslin-Voevodsky's pullbacks of relative cycles coincide with the Koll\'ar's ones in our case. Thus, for each $S$, one has the bijection
  $$
  \Upsilon ^p_d(X,S):\bcC ^p_d(X\times S/S)
  \stackrel{\sim }{\lra }
  \Hom (S,C^p_d(X))\; ,
  $$
functorial in $S$. Moreover, these bijections are also functorial in $X$ by Corollary 3.6.3 in \cite{SuslinVoevodskyChow}.

Notice that if $d=0$, then, by convention, $C^p_0(X)=\Spec (k)$ and the unique $k$-point of $\Spec (k)$ is the neutral element $0$ of the free monoid $\bcC ^p(X)$. Since $0$ can be also considered as the empty codimension $p$ relative cycle over $\Spec (k)$, we in fact identify the unique $k$-point of $\Spec (k)$ with the empty relative cycle.

It is trivial but worth noticing that if $k'$ is another field and $\alpha :k\stackrel{\sim }{\to }k'$ is an isomorphism of fields, the scheme $C_r(X')$ is the pull-back of the scheme $C_r(X)$ with respect to the morphism $\Spec (\alpha )$, where $X'$ is the pull-back of $X$, and the corresponding morphism from $C^p(X')$ to $C^p(X)$ is an isomorphism of schemes. The bijections $\Upsilon ^p(X',S)$ and $\Upsilon ^p(X,S)$ commute through the obvious isomorphisms on monoids and Hom-sets, induced by the isomorphism $\alpha $.

Next, for a commutative (additive) monoid $M$, its completion $M^+$ can be constructed as the quotient of $M\oplus M$ by the image of the diagonal embedding. Let
  $$
  \tau :M\oplus M\to M^+
  $$
be the corresponding quotient homomorphism, and let
  $$
  \nu :M\to M^+
  $$
be the composition of the embedding of $M$ as one of the two direct summands and the homomorphism $\tau $. Then $\nu $ possesses the obvious universal property, and for any $(a,b)$ in $M\oplus M$ the value $\tau (a,b)$ is the difference $\nu (a)-\nu (b)$. If $M$ is a cancellation monoid then $\nu $ is injective, and we can identify $M$ with its image in $M^+$. Modulo this identification,
  $$
  \tau (a,b)=a-b\; .
  $$
In particular, we can consider the presheaf $\bcZ ^p(X)$ of abelian groups on $\Noe $, such that for each $S$ the group of sections
  $$
  \bcZ ^p(X\times S/S)=\bcC ^p(X\times S/S)^+
  $$
is the group completion of the monoid $\bcC ^p(X\times S/S)$.

Identifying schemes with representable presheaves, we identify the presheaf $\bcC ^p_d(X)$ with the Chow scheme $C^p_d(X)$. Looking at $\bcC ^p(X)=\coprod _d\bcC ^p_d(X)$ as a monoid object in the category of presheaves, we may also interpret the presheaf $\bcZ ^p(X)$ as the group completion $\bcC ^p(X)^+$ of this monoid in the category of presheaves.

Let $S$ and $Y$ be two Noetherian schemes over $k$. Define a functor $\cHom (S,Y)$ on $\Noe $ sending a scheme $T$ to the set $\cHom (S,Y)(T)$ of morphisms from $S\times T$ to $Y\times T$ over $T$. The graphs of such morphisms give us an embedding of $\cHom (S,Y)$ into the Hilbert functor $\cHilb (S\times Y)$. If $S$ and $Y$ are projective over $k$, the latter functor is representable by the projective Hilbert scheme $\Hilb (S\times Y)$ and $\cHom (S,Y)$ is representable by an open subscheme $\Hom (S,Y)$ in $\Hilb (S\times Y)$. Given a positive integer $d$, one can choose a suitable Hilbert polynomial $\Phi _d$, such that the intersection $\cHom ^d(S,Y)$ of $\cHom (S,Y)$ with $\cHilb _{\Phi _d}(S\times Y)$ inside $\cHilb (S\times Y)$ is representable by an open subscheme $\Hom ^d(S,Y)$ in the projective scheme $\Hilb _{\Phi _d}(S\times Y)$, and if $S=\PR ^1$ then $\Hom ^d(\PR ^1,Y)$ is a quasi-projective scheme over $k$ parametrizing rational curves of degree $d$ in $Y$. By the universal property of fibred products over $k$ one has the natural bijection between $\cHom (S,Y)(T)$ and $\Hom (T\times S,Y)$. This gives us the adjunction
  $$
  \Hom (T\times S,Y)\simeq \Hom (T,\Hom (S,Y))
  $$
and the corresponding regular evaluation morphism from the scheme $\Hom (S,Y)\times S$ to the scheme $Y$. The latter morphism induces the regular evaluation morphism of quasi-projective schemes
  $$
  e_{S,Y}:\Hom ^d(S,Y)\times S\to Y\; ,
  $$
for each positive integer $d$. In particular, if $P$ is a closed point of $\PR ^1$ one has the regular evaluation morphism
  $$
  e_P:\Hom ^d(\PR ^1,Y)\to Y\; ,
  $$
sending $f$ to $f(P)$. This all can be found in the first chapter of \cite{Kollar}.

Let now $X$ be a smooth projective variety embedded into $\PR ^n$ over $k$. Let $A$ and $A'$ be two algebraic cycles of codimension $p$ on $X$. The cycle $A$ is rationally equivalent to the cycle $A'$ if and only if there exists an effective codimension $p$ algebraic cycle $Z$ on $X\times \PR ^1$ and an effective codimension $p$ algebraic cycle $B$ on $X$, such that $$
  Z(0)=A+B
  \qqand
  Z(\infty )=A'+B\; .
  $$
Assume that $A$ is rationally equivalent to $A'$, and let
  $$
  f_Z=\Upsilon (Z)
  $$
and
  $$
  f_{B\times \PR ^1}=\Upsilon (B\times \PR ^1)
  $$
be two regular morphisms from $\PR ^1$ to $\bcC ^p(X)$, where
  $$
  \Upsilon =\Upsilon ^p(X,\PR ^1):
  \bcC ^p(X\times \PR ^1/\PR ^1)\to
  \Hom (\PR ^1,C^p(X))
  $$
is the functorial bijection considered above. Let also
  $$
  f=f_Z\oplus f_{B\times \PR ^1}:
  \PR ^1\to C^p(X)\oplus C^p(X)
  $$
be the morphism generated by $f_Z$ and $f_{B\times \PR ^1}$. Since $C^p(X)$ is a cancellation monoid, for any two elements $a,b\in C^p(X)$ the value $\tau (a,b)$ in $C^p(X)^+$ is $a-b$, after the identification of $C^p(X)$ with its image in $C^p(X)^+$ under the injective homomorphism $\nu $ from $C^p(X)$ to $C^p(X)^+$. Then
  $$
  \tau f(0)=\tau (f_Z(0),f_{B\times \PR ^1})=
  f_Z(0)-f_{B\times \PR ^1}(0)=Z(0)-B=A
  $$
and
  $$
  \tau f(\infty )=\tau (f_Z(\infty ),f_{B\times \PR ^1})=
  f_Z(0)-f_{B\times \PR ^1}(\infty )=Z(\infty )-B=A'\; .
  $$

Vice versa, suppose there is a regular morphism $f=f_1\oplus f_2$ from $\PR ^1$ to the direct sum $C^p(X)\oplus C^p(X)$, such that $\tau f(0)=A$ and $\tau f(\infty )=A'$. Let $Z_1$ and $Z_2$ be two algebraic cycles in $\bcC ^p(X\times \PR ^1/\PR ^1)$, such that $\Upsilon (Z_i)=f_i$ for $i=1,2$, and let $Z=Z_1-Z_2$. Then $Z(0)=A$ and $Z(\infty )=A'$, which means that $A$ is rationally equivalent to $A'$.

For any non-negative integers $d_1,\dots ,d_s$ let
  $$
  C^p_{d_1,\dots ,d_s}(X)=
  C^p_{d_1}(X)\times \dots \times C^p_{d_s}(X)
  $$
be the fibred product over the ground field $k$. For any degree $d\geq 0$ let
  $$
  W_d=
  \{ (A,B)\in C^p_{d,d}(X)\; |\; A\stackrel{\rat }{\sim }B\}
  $$
be the Zariski closed subset in $C^p_{d,d}(X)$ determined by ordered pairs $(A,B)$ of closed points in $C^p_d(X)$, such that the cycle $A$ is rationally equivalent to the cycle $B$ on $X$. For any non-negative $u$ and positive $v$ let also
  $$
  W_d^{u,v}=
  $$
  $$
  \{ (A,B)\in C^p_{d,d}(X)\; |\; \exists
  f\in \Hom ^v(\PR ^1,C^p_{d+u,u}(X)),\; \hbox{s.th.}\;
  \tau f(0)=A, \tau f(\infty )=B\}
  $$
Then
  $$
  W_d^{u,v}\subset W_d
  \qqand
  W_d=\cup _{u,v}W_d^{u,v}\; .
  $$
Let also $\bar W_d^{u,v}$ be the Zariski closure of the set $W_d^{u,v}$ in the scheme $C^p_{d,d}(X)$.

\begin{proposition}
\label{l1}
For any $d$, $u$ and $v$, the set $W_d^{u,v}$ is a quasi-projective subscheme in $C^p_{d,d}(X)$ whose Zariski closure $\bar W_d^{u,v}$ is contained in $W_d$.
\end{proposition}

\begin{pf} To prove the proposition all we need is to extend the arguments in \cite{Roitman} from zero-cycles and symmetric powers to codimension $p$ cycles and Chow varieties.

Let
  $$
  e:\Hom ^v(\PR ^1,C^p_{d+u,u}(X))\to
  C^p_{d+u,u,d+u,u}(X)
  $$
be the evaluation morphism sending $f$ to the ordered pair $(f(0),f(\infty ))$, and let
  $$
  s:C^p_{d,u,d,u}(X)\to C^p_{d+u,u,d+u,u}(X)
  $$
be the regular morphism given by the formula
  $$
  s(A,C,B,D)=(A+C,C,B+D,D)\; .
  $$
The two morphisms $e$ and $s$ allow us to take the product
  $$
  V=\Hom ^v(\PR ^1,C^p_{d+u,u}(X))
  \times _{C^p_{d+u,u,d+u,u}(X)}C^p_{d,u,d,u}(X)\; ,
  $$
which is a closed subvariety in the product
  $$
  \Hom ^v(\PR ^1,C ^p_{d+u,u}(X))\times C^p_{d,u,d,u}(X)
  $$
over $\Spec (k)$ consisting of quintuples $(f,A,C,B,D)$ with $e(f)=s(A,C,B,D)$, i.e.
  $$
  (f(0),f(\infty ))=(A+C,C,B+D,D)\; .
  $$
The latter equality gives that
  $$
  \pr _{2,3}(V)\subset W_d^{u,v}\; ,
  $$
where $\pr _{2,3}$ is the projection of $\Hom ^v(\PR ^1,C^p_{d+u,u}(X))\times C^p_{d,u,d,u}(X)$ onto $C^p_{d,d}(X)$.

Vice versa, if $(A,B)$ is a closed point of $W_d^{u,v}$, there exists a regular morphism
  $$
  f\in \Hom ^v(\PR ^1,C^p_{d+u,u}(X))
  $$
with
  $$
  \tau f(0)=A\qqand \tau f(\infty )=B\; .
  $$
Let $f(0)=(C',C)$ and $f(\infty )=(D',D)$. Then
  $$
  \tau f(0)=C'-C=A
  $$
and
  $$
  \tau f(\infty )=D'-D=B
  $$
in the completed monoid $\bcZ ^p(X)=C^p(X)^+$. This means that there exist effective codimension $p$ algebraic cycles $M$ and $N$ on $X$, such that
  $$
  C'+M=C+A+M
  \qqand
  D'+N=D+B+N
  $$
in $C^p(X)$. Since $C^p(X)$ is a free monoid, it possesses the cancellation property. Therefore, the latter two equalities imply the equalities
  $$
  C'=C+A
  \qqand
  D'=D+B
  $$
respectively. This yields
  $$
  e(f)=s(A,C,B,D)\; ,
  $$
whence
  $$
  (f,A,C,B,D)\in V\; .
  $$
It means that $(A,B)$ is in $\pr _{2,3}(V)$.

Thus,
  $$
  \pr _{2,3}(V)=W_d^{u,v}\; .
  $$
Being the image of a quasi-projective variety under the projection $\pr _{2,3}$ the set $W_d^{u,v}$ is itself a quasi-projective variety.

Let
  $$
  \tilde s:C^p_{d,d,u,u}(X)\to C^p_{d+u,d+u,u,u}(X)
  $$
be the morphism obtained by composing and precomposing $s$ with the transposition of the second and third factors in the domain and codomain of $s$. Then
  $$
  W_d=\pr _{1,2}(\tilde s^{-1}(W_{d+u}\times W_u))\; .
  $$
Let $(A,B,C,D)$ be a closed point in $C^p_{d,d,u,u}(X)$, such that the value
  $$
  \tilde s(A,B,C,D)=(A+C,B+D,C,D)
  $$
is in $W_{d+u}^{0,v}\times W_u^{0,v}$. It means that there exist two regular morphisms
  $$
  g\in \Hom ^v(\PR ^1,C^p_{d+u}(X))
  \qqand
  h\in \Hom ^v(\PR ^1,C^p_u(X))
  $$
with
  $$
  g(0)=A+C\; ,\; \; g(\infty )=B+D\; ,\; \; h(0)=C
  \qand h(\infty )=D\; .
  $$
Then
  $$
  f=g\times h\in \Hom ^v(\PR ^1,C^p_{d+u,u}(X))\; ,
  $$
  $$
  f(0)=(A+C,C)\qqand f(\infty )=(B+D,D)\; .
  $$
Hence, $\tau f(0)=A$ and $\tau f(\infty )=B$. It means that $(A,B)\in W_d^{u,v}$. We have shown that
  $$
  \pr _{1,2}(\tilde s^{-1}(W_{d+u}^{0,v}\times W_u^{0,v}))\subset W_d^{u,v}\; .
  $$

Vice versa, suppose $(A,B)\in W_d^{u,v}$, and let $f$ be a morphism from $\PR ^1$ to $C^p_{d+u,u}(X)$ with $\tau f(0)=A$ and $\tau f(\infty )=B$. Composing $f$ with the projections of $C^p_{d+u,u}(X)$ onto $C^p_{d+u}(X)$ and $C^p_u(X)$ one can easily show that
  $$
  W_d^{u,v}\subset
  \pr _{1,2}(\tilde s^{-1}(W_{d+u}^{0,v}\times W_u^{0,v}))\; .
  $$

Thus,
  $$
  W_d^{u,v}=
  \pr _{1,2}(\tilde s^{-1}(W_{d+u}^{0,v}\times W_u^{0,v}))\; .
  $$
Since $\tilde s$ is continuous and $\pr _{1,2}$ is proper, we obtain that
  $$
  \bar W_d^{u,v}=
  \pr _{1,2}(\tilde s^{-1}(\bar W_{d+u}^{0,v}
  \times \bar W_u^{0,v}))\; .
   $$
This gives us that in order to prove the second assertion of the proposition it suffices to show that $\bar W_d^{0,v}$ is contained in $W_d$.

Let $(A,B)$ be a closed point of $\bar W_d^{0,v}$. If $(A,B)$ is in $W_d^{0,v}$, then it is also in $W_d$. Suppose that $(A,B)$ is in $\bar W_d^{0,v}\smallsetminus W_d^{0,v}$. Let $W$ be an irreducible component of the quasi-projective variety $W_d^{0,v}$ whose Zariski closure $\bar W$ contains the point $(A,B)$. Let $U$ be an affine neighbourhood of $(A,B)$ in $\bar W$. Since $(A,B)$ is in the closure of $W$ the set $U\cap W$ is non-empty. Let $C$ be an irreducible curve passing through $(A,B)$ in $U$ and let $\bar C$ be the Zariski closure of $C$ in $\bar W$. The evaluation regular morphisms
  $$
  e_0:\Hom ^v(\PR ^1,C_d^p(X))\to C_d^p(X)
  \qand
  e_{\infty }:\Hom ^v(\PR ^1,C_d^p(X))\to C_d^p(X)
  $$
give the regular morphism
  $$
  e_{0,\infty }:\Hom ^v(\PR ^1,C_d^p(X))\to C_{d,d}^p(X)\; .
  $$
Then $W_d^{0,v}$ is exactly the image of the regular morphism $e_{0,\infty }$, and we can choose a quasi-projective curve $T$ in $\Hom ^v(\PR ^1,C_d^p(X))$, such that the closure of the image $e_{0,\infty }(T)$ is $\bar C$. Since $\Hom ^v(\PR ^1,\bcC _d^p(X))$ is a quasi-projective variety, we can embed it into some projective space $\PR ^m$. Let $\bar T$ be the closure of $T$ in $\PR ^m$, let $\tilde T$ be the normalization of $\bar T$ and let $\tilde T_0$ be the pre-image of $T$ in $\tilde T$. Consider the composition
  $$
  f_0:\tilde T_0\times \PR ^1\to T\times \PR ^1\subset
  \Hom ^v(\PR ^1,C_d^p(X))\times \PR ^1
  \stackrel{e}{\lra }C_d^p(X)\; ,
  $$
where $e$ is the evaluation morphism $e_{\PR ^1,C_d^p(X)}$. The regular morphism $f_0$ defines a rational map
  $$
  f:\tilde T\times \PR ^1\dasharrow C_d^p(X)
  $$
Since $\tilde T$ is a smooth projective curve, the product $\tilde T\times \PR ^1$ is a smooth projective surface over the ground field. Under this condition there exists a finite chain of $\sigma $-processes $(\tilde T\times \PR ^1)'\to \tilde T\times \PR ^1$ resolving indeterminacy of $f$ and giving a regular morphism
  $$
  f':(\tilde T\times \PR ^1)'\to C_d^p(X)\; .
  $$
The regular morphism $\tilde T_0\to T\to \bar C$ extends to the regular morphism $\tilde T\to \bar C$. Let $P$ be a point in the fibre of this morphism at $(A,B)$. For any closed point $Q$ on $\PR ^1$ the restriction $f|_{\tilde T\times \{ Q\} }$ of the rational map $f$ onto $\tilde T\times \{ Q\} \simeq \tilde T$ is regular on the whole curve $\tilde T$, because $\tilde T$ is smooth. Then
  $$
  (f|_{\tilde T\times \{ 0\} })(P)=A
  $$
and
  $$
  (f|_{\tilde T\times \{ \infty \} })(P)=B\; .
  $$
It means that the points $A$ and $B$ are connected by a finite collection of curves which are the images of rational curves on $(\tilde T\times \PR ^1)'$ under the regular morphism $f'$. In turn, it follows that $A$ is rationally equivalent to $B$, whence
  $$
  (A,B)\in W_d\; .
  $$
\end{pf}

In what follows, for any equi-dimensional algebraic scheme $V$ over $k$, let $CH^p(V)$ be the Chow group, with coefficients in $\ZZ $, of codimension $p$ algebraic cycles modulo rational equivalence on $V$. If a closed embedding $V\subset \PR ^m$ is fixed, let $CH^p(V)_{\deg =0}$ be the subgroup generated by cycles classes of degree $0$ in $CH^p(V)$. Then, for any two nonnegative integer $d$, we have a map
  $$
  \theta ^p_d:C^p_{d,d}(X)\to CH^p(X)_{\deg =0}
  $$
sending an ordered pair $(A,B)$ of closed points on the Chow variety $C^p_d(X)$ to the class of the difference $Z_A-Z_B$, where $Z_A$ and $Z_B$ are codimension $p$ algebraic cycles on $X$ corresponding to the points $A$ and $B$ respectively.

\begin{corollary}
\label{maincor}
$(\theta ^p_d)^{-1}(0)$ is a countable union of irreducible Zariski closed subsets in the Chow scheme $C^p_{d,d}(X)$.
\end{corollary}

\begin{pf}
Proposition \ref{l1} gives that $W_d$ is the countable union of Zariski closed sets $\bar W_d^{u,v}$ over $u$ and $v$. This completes the proof.
\end{pf}

\section{Countability lemmas and the abelian variety $A_0$}
\label{proofofA}

The purpose of this section is to prove Theorem A in Introduction, which generalizes the argument from Section 10.1.2 in the second volume of the book \cite{Voisin: The Book} (see also pp 304 - 305 in loc.cit.). Theorem A introduces the abelian variety $A_0$ which is of key importance to the whole approach. We shall also prove that $A_0$ is a subvariety of the abelian variety $A_1$ introduced in Section \ref{setting&assumptions}. Since now and throughout the rest of the paper we shall assume that the ground field $k$ is uncountable.

\begin{lemma}
\label{zayac}
Let $V$ be an irreducible quasi-projective algebraic variety over $k$. Then $V$ cannot be written as a countable union of its Zariski closed subsets, each of which is not the whole $V$.
\end{lemma}

\begin{pf}
Since $V$ is supposed to be irreducible, without loss of generality we may assume that $V$ is affine. Let $d$ be the dimension of $V$ and suppose $V=\cup _{n\in \NN }V_n$ is the union of closed subsets $V_n$ in $V$, such that $V_n\neq V$ for each $n$. By Emmy Noether's lemma, there exists a finite surjective morphism $f:V\to \AF ^d$ over $k$. Let $W_n$ be the image of $V_n$ under $f$. Since $f$ is finite, it is proper. Therefore, $W_n$ are closed in $\AF ^d$, so that we obtain that the affine space $\AF ^d$ is the union of $W_n$'s. Since the ground field is uncountable, the set of all hyperplanes in $\AF ^d$ is uncountable. Therefore, there exists a hyperplane $H$, such that $W_n\not \subset H$ for any index $n$. Induction reduces the assertion of the lemma to the case when $d=1$.
\end{pf}

A countable union $V=\cup _{n\in \NN }V_n$ of algebraic varieties will be called irredundant if $V_n$ is irreducible for each $n$ and $V_m\not \subset V_n$ for $m\neq n$. In an irredundant decomposition, the sets $V_n$ will be called $c$-components of $V$.

\begin{lemma}
\label{belka}
Let $V$ be a countable union of algebraic varieties over an uncountable algebraically closed ground field. Then $V$ admits an irredundant decomposition, and such an irredundant decomposition is unique.
\end{lemma}

\begin{pf}
Let $V=\cup _{n\in \NN }V'_n$ be a countable union of algebraic varieties over $k$. For each $n$ let $V'_n=V'_{n,1}\cup \dots \cup V'_{n,r_n}$ be the irreducible components of $V'_n$. Ignoring all components $V'_{m,i}$ with  $V'_{m,i}\subset V'_{n,j}$ for some $n$ and $j$ we end up with a irredundant decomposition. Having two irredundant decompositions $V=\cup _{n\in \NN }V_n$ and $V=\cup _{n\in \NN }W_n$, suppose there exists $V_m$ such that $V_m$ is not contained in $W_n$ for any $n$. Then $V_m$ is the union of the closed subsets $V_m\cap W_n$, each of which is not $V_m$. This contradicts to Lemma \ref{zayac}. Therefore, any $V_m$ is contained in some $W_n$. By symmetry, any $W_n$ is in $V_l$ for some $l$. Then $V_m\subset V_l$. By irredundancy, $l=m$ and $V_m=W_n$.
\end{pf}

\begin{lemma}
\label{lisa}
Let $A$ be an abelian variety over $k$, and let $K$ be a subgroup which can be represented as a countable union of Zariski closed subsets in $A$. Then the irredundant decomposition of $K$ contains a unique irreducible component passing through $0$, and this component is an abelian subvariety in $A$.
\end{lemma}

\begin{pf}
Let $K=\cup _{n\in \NN }K_n$ be the irredundant decomposition of $K$, which exists by Lemma \ref{belka}. Since $0\in K$, there exists at least one component in the irredundant decomposition, which contains $0$. Suppose there are $s$ components $K_1,\dots ,K_s$ containing $0$ and $s>1$. The summation in $K$ gives the regular morphism from the product $K_1\times \dots \times K_s$ into $A$, whose image is the irreducible Zariski closed subset $K_1+\dots +K_s$ in $A$. By Lemma \ref{zayac}, there exists $n\in \NN $, such that $K_1+\dots +K_s\subset K_n$ and so $K_i\subset K_n$ for each $1\leq i\leq n$. By irredundancy, $s=1$, which contradicts to the assumption $s>1$.

After renumbering of the components, we may assume that $0\in K_0$. If $K_0=\{ 0\} $, then $K_0$ is trivially an abelian variety. Suppose $K_0\neq \{ 0\} $ and take a non-trivial element $x$ in $K_0$. Since $-x+K_0$ is irreducible, it must be in some $K_n$ by Lemma \ref{zayac}. As $0\in -x+K_0$ it follows that $0\in K_n$ and so $n=0$. It follows that $-x\in K_0$. Similarly, since $K_0+K_0$ is irreducible and contains $K_0$, we see that $K_0+K_0=K_0$. Being a Zariski closed abelian subgroup in $A$, the set $K_0$ is an abelian subvariety in $A$.
\end{pf}


Now we are ready to prove Theorem A in Introduction.

\begin{theorem}
\label{kot}
In terms above, there exists an abelian subvariety $A_0$ and a countable subset of closed points in $A$, such that
  $$
  K=\cup _{x\in \Xi }(x+A_0)
  $$
inside the abelian variety $A$.
\end{theorem}

\begin{pf}
Let $n$ be the dimension of $Y$. Since $\bcC _{n-p}$ is a functor on Noetherian schemes over $k$, the closed immersion $r$ induces a morphism from $\bcC _{n-p}(Y)$ to $\bcC _{n-p}(X)$ in the category presheaves on seminormal schemes. In upper indices, this is a morphism from $\bcC ^p(Y)$ to $\bcC ^{p+e}(X)$. Passing to Chow schemes, we obtain a regular morphism
  $$
  r_*:C^p(Y)\to C^{p+e}(X)
  $$
of projective schemes over $k$. Since $X$ is embedded into $\PR ^m$, and since $Y$ is a closed subvariety in $X$, the morphism $r_*$ induces the morphisms
  $$
  r_*:C^p_d(Y)\to C^{p+e}_d(X)\; ,
  $$
one for each degree $d$. Taking in to account Assumption 2, we obtain the obvious commutative diagram
  $$
  \diagram
  \coprod _dC^p_{d,d}(Y)
  \ar[dd]_-{\coprod \theta ^p_d}
  \ar[rr]^-{\coprod _dr_*} & &
  \coprod _dC^{p+e}_{d,d}(X)
  \ar[dd]^-{\coprod \theta ^{p+e}_d} \\ \\
  A^p(Y) \ar[rr]^-{r_*} & & CH^{p+e}(X)_{\deg =0}
  \enddiagram
  $$

Since $A^p(Y)$ is weakly representable, there exists a smooth projective curve $\Gamma $ over $k$, and an algebraic cycle $Z$ of codimension $p$ on $\Gamma \times Y$, such that the induced homomorphism
  $$
  Z_*:A^1(\Gamma )\to A^p(Y)
  $$
is surjective. On the other hand, since $A^1(\Gamma )$ is representable by the Jacobian of $\Gamma $, the map
  $$
  \theta ^1_d:C^1_{d,d}(\Gamma )\to A^1(\Gamma )
  $$
is surjective for big enough $d$. Using these two facts one can show that the right vertical arrow of the above commutative diagram is surjective. Then the kernel $K$ of the bottom horizontal homomorphism is the image under the map $\coprod \theta ^p_d$ of the preimage of $0$ under the composition of the maps $\coprod _dr_*$ and $\coprod \theta ^{p+e}_d$. Corollary \ref{maincor} implies that the latter preimage is the coproduct of countable unions of Zariski closed subsets in the schemes $C^p_{d,d}(Y)$.

Now, consider the composition
  $$
  \psi \circ \theta ^p_d:C^p_{d,d}(Y)\to A^p(Y)\to A\; ,
  $$
for each number $d$. By the definition of the regularity of $\psi $, this composition is a regular morphism of schemes. Since these schemes are projective, the composition is proper. It follows that the subgroup $\psi (K)$ is a countable union of Zariski closed subsets in the abelian variety $A$.

For simplicity of notation, identify $\psi (K)$ and $K$. By Lemma \ref{belka}, the set $K$ admits a unique irredundant decomposition inside the abelian variety $A$, say
  $$
  K=\cup _{n\in \NN }K_n\; .
  $$
Let $A_0$ be the unique component of that decomposition passing through $0$, which is an abelian subvariety in $A$ by Lemma \ref{lisa}. For any $x$ in $K$ the set $x+A_0$ is an irreducible Zariski closed subset in $K$. Since $K$ coincides with $\cup _{x\in K}(x+A_0)$, ignoring each set $x+A_0$ which is a subsets in $y+A_0$ for some $y\in K$, we can find a subset $\Xi $ in $K$, such that
  $$
  K=\cup _{x\in \Xi }(x+A_0)\; ,
  $$
and for any two elements $x,y\in \Xi $ the irreducible sets $x+A_0$ and $y+A_0$ are not contained one in another. Since $x+A_0$ is irreducible, it is contained in $K_n$ for some $n$ by Lemma \ref{zayac}. Then $A_0\subset -x+K_n$. Similarly, $-x+K_n\subset K_l$ for some $l$, so that $K_l=A_0$ by Lemma \ref{lisa}. This yields
  $$
  x+A_0=K_n\; ,
  $$
for each $x\in \Xi $. It means that the set $\Xi $ is countable.
\end{pf}

Let us also prove that the abelian variety $A_0$, provided by Theorem \ref{kot}, is contained in the abelian variety $A_1$ introduced in Section \ref{setting&assumptions}. Choose an ample line bundle $L$ on the abelian variety $A$. Let
  $$
  i:A_0\to A
  $$
be the closed embedding of $A_0$ into $A$, and let $L_0$ be the pull-back of $L$ to $A_0$ under the embedding $i$. Define the homomorphism $\zeta $ on divisors via the commutative diagram
  $$
  \diagram
  A^1(A_0)
  \ar[dd]_-{(\lambda _{L_0})_*}
  \ar[rr]^-{\zeta } & &
  A^1(A) \\ \\
  A^1(A_0^{\vee })
  \ar[rr]^-{{i^{\vee }}^*}  & &
  A^1(A^{\vee }) \ar[uu]^-{\lambda _L^*}
  \enddiagram
  $$
Similarly, we define the homomorphism $\zeta _{\ZZ _l}$ on cohomology by means of the commutative diagram
  $$
  \diagram
  H^1_{\et }(A_0,\ZZ _l)
  \ar[rr]^-{\zeta _{\ZZ _l}}
  \ar[dd]^-{{\lambda _{L_0}}_*} & &
  H^1_{\et }(A,\ZZ _l) \\ \\
  H^1_{\et }(A_0^{\vee },\ZZ _l)
  \ar[rr]^-{{i^{\vee }}^*} & &
  H^1_{\et }(A^{\vee },\ZZ _l)
  \ar[uu]_-{\lambda _L^*}
  \enddiagram
  $$
and analogously for the homomorphism
   $$
   \zeta _{\QQ _l/\ZZ _l}:
   H^1_{\et }(A_0,\QQ _l/\ZZ _l)\to
   H^1_{\et }(A_0,\QQ _l/\ZZ _l)\; .
   $$
The homomorphism $\zeta _{\ZZ _l}$ induces the injective  homomorphisms
   $$
   \zeta _{\QQ _l}:
   H^1_{\et }(A_0,\QQ _l)=
   H^1_{\et }(A_0,\ZZ _l)\otimes \QQ _l
   \to H^1_{\et }(A_0,\QQ _l)=
   H^1_{\et }(A,\ZZ _l)\otimes \QQ _l
   $$
and
   $$
   \zeta _{\ZZ _l}\otimes \QQ _l/\ZZ _l:
   H^1_{\et }(A_0,\ZZ _l)\otimes \QQ _l/\ZZ _l\to
   H^1_{\et }(A_0,\ZZ _l)\otimes \QQ _l/\ZZ _l\; .
   $$

\begin{proposition}
\label{kutki}
The image of the composition
  $$
  H^1_{\et }(A_0,\QQ _l(1-p))
  \stackrel{\zeta _{\QQ _l}}{\lra }
  H^1_{\et }(A,\QQ _l(1-p))\stackrel{w_*}{\lra }
  H^{2p-1}_{\et }(Y,\QQ _l)
  $$
is contained in the kernel of the Gysin homomorphism
  $$
  H^{2p-1}_{\et }(Y,\QQ _l)
  \stackrel{r_*}{\lra }
  H^{2(p+e)-1}_{\et }(X,\QQ _l)\; .
  $$
\end{proposition}

\begin{pf}
For the proof we will be using Bloch's $l$-adic Abel-Jacobi maps. For any abelian group $A$, a prime $l$ and positive integer $n$ let $A_{l^n}$ be the kernel of the multiplication by $l^n$ endomorphism of $A$ and let $A(l)$ be the $l$-primary part of $A$, i.e. the union of the groups $A_{l^n}$ for all $n$. For any smooth projective variety $V$ over $k$, there is a canonical homomorphism
  $$
  \lambda ^p_l(V): CH^p(V)(l)\to
  H^{2p-1}_{\et }(V,\QQ _l/\ZZ _l(p))\; ,
  $$
constructed by Bloch in \cite{Bloch79}. The homomorphisms $\lambda^p_l(V)$ are functorial with respect to the action of correspondences between smooth projective varieties over $k$. Moreover, the homomorphisms
  $$
  \lambda ^1_l(V):CH^1(V)(l)\to H^1_{\et }(V,\QQ _l/\ZZ _l(1))
  $$
are all isomorphisms, loc.cit.

Since $A_0$ and $A$ are abelian varieties, their N\'eron-Severi groups are torsion free. It follows that $CH^1(A_0)(l)=A^1(A_0)(l)$ and $CH^1(A)(l)=A^1(A)(l)$, so that we actually have the isomorphism $\lambda ^1_l(A_0)$ between $A^1(A_0)(l)$ and $H^1_{\et }(A_0,\QQ _l/\ZZ _l(1))$, and the isomorphism $\lambda ^1_l(A)$ between $A^1(A)(l)$ and $H^1_{\et }(A,\QQ _l/\ZZ _l(1))$. Similarly, one has the isomorphisms $\lambda ^1_l(A_0^{\vee })$ and $\lambda ^1_l(A^{\vee })$ for the dual abelian varieties.

The functorial properties of Bloch's maps $\lambda ^1_l$ give us that the diagram
  \begin{equation}
  \label{candibobior}
  \diagram
  A^1(A_0)(l) \ar[dd]^{\sim }_-{\lambda ^1_l(A_0)}
  \ar[rr]^-{\zeta } & &
  A^1(A)(l) \ar[dd]_-{\lambda ^1_l(A)}^-{\sim }
  \\ \\
  H^1_{\et }(A_0,\QQ _l/\ZZ _l(1))
  \ar[rr]^-{\zeta _{\QQ _l/\ZZ _l}} & &
  H^1_{\et }(A,\QQ _l/\ZZ _l(1))
  \enddiagram
  \end{equation}
is commutative.

For a smooth projective $V$ over $k$, one has the homomorphisms
  $$
  \varrho ^{i,j}_l(V):H^i_{\et }(V,\ZZ _l(j))\otimes \QQ _l/\ZZ _l\to H^i_{\et }(V,\QQ _l/\ZZ _l(j))\; ,
  $$
with finite kernels and cokernels, considered in \cite{GorchGulFano}. In particular, we have the commutative diagram
  \begin{equation}
  \label{rhomorphism}
  \diagram
  H^1_{\et }(A_0,\QQ _l/\ZZ _l(1))
  \ar[rr]^-{\zeta _{\QQ _l/\ZZ _l}} & &
  H^1_{\et }(A,\QQ _l/\ZZ _l(1)) \\ \\
  H^1_{\et }(A_0,\ZZ _l(1))\otimes \QQ _l/\ZZ _l \ar[uu]_-{\varrho _l^{1,1}(A_0)} \ar[rr]^-{\zeta _{\ZZ _l}
  \otimes \QQ _l/\ZZ _l} & &
  H^1_{\et }(A,\ZZ _l(1))\otimes \QQ _l/\ZZ _l
  \ar[uu]_-{\varrho _l^{1,1}(A)}
  \enddiagram
  \end{equation}

Let $\sigma :A_0\stackrel{\sim }{\to }A^1(A_0^{\vee })$ and $\sigma :A\stackrel{\sim }{\to }A^1(A^{\vee })$ be the autoduality isomorphisms. The morphism of motives
  $$
  w:M(A)\otimes \Le ^{p-1}\to M(Y)
  $$
induces the homomorphism $w_*:A^1(A)\to A^p(Y)$ on Chow groups. A straightforward verification shows that the diagram
  \begin{equation}
  \label{iozhik}
  \diagram
  A^1(A)(l) \ar[rr]^-{w_*} & & A^p(Y)(l) \\ \\
  A^1(A^{\vee })(l) \ar[uu]_-{\lambda _L^*} & &
  A(l) \ar[ll]_-{\sigma } \ar[uu]_-{(\psi ^p_Y)^{-1}}
  \enddiagram
  \end{equation}
is commutative.

The homomorphism $w_*:A^1(A)\to A^p(Y)$ on Chow groups and the homomorphism $w_*:H^1_{\et }(A,\QQ _l/\ZZ _l(1))\to H^{2p-1}_{\et }(Y,\QQ _l/\ZZ _l(p))$ induced by $w$ on cohomology fit into the commutative diagram
  \begin{equation}
  \label{mohnatyiles}
  \diagram
  A^1(A)(l) \ar[dd]^{\sim }_-{\lambda ^1_l(A)} \ar[rr]^-{w_*} & &
  A^p(Y)(l) \ar[dd]_-{\lambda ^p_l(Y)}^-{} \\ \\
  H^1_{\et }(A,\QQ _l/\ZZ _l(1)) \ar[rr]^-{w_*} & &
  H^{2p-1}_{\et }(Y,\QQ _l/\ZZ _l(p))
  \enddiagram
  \end{equation}

The commutativity of the diagrams (\ref{candibobior}), (\ref{rhomorphism}), (\ref{iozhik}) and (\ref{mohnatyiles}), the definition of the abelian variety $A_0$ and easy diagram chase over the obvious commutative diagrams with Gysin mappings show that the image of the triple composition
  $$
  r_*\circ w_*\circ (\zeta _{\ZZ _l}\otimes \QQ _l/\ZZ _l):
  H^1_{\et }(A_0,\ZZ _l(1))\otimes \QQ _l/\ZZ _l\to H^{2(p+e)-1}_{\et }(X,\ZZ _l(p+e))\otimes \QQ _l/\ZZ _l
  $$
is contained in the kernel of the homomorphism
  $$
  \varrho ^{2(p+e)-1,p+e}_l(X):
  H^{2(p+e)-1}_{\et }(X,\ZZ _l(p+e))\otimes \QQ _l/\ZZ _l\to H^{2(p+e)-1}_{\et }(X,\QQ _l/\ZZ _l(p+e))\; .
  $$
Since the latter kernel is finite, the image of the composition $r_*\circ w_*\circ (\zeta _{\ZZ _l}\otimes \QQ _l/\ZZ _l)$ is finite too. Since the \'etale cohomology groups of smooth projective varieties with $\ZZ _l$-coefficients are finitely generated $\ZZ _l$-modules, it follows that the image of the triple composition
  $$
  r_*\circ w_*\circ \zeta _{\QQ _l}:
  H^1_{\et }(A_0,\QQ _l(1))\to
  H^{2(p+e)-1}_{\et }(X,\QQ _l(p+e))
  $$
is zero, which finishes the proof of the proposition.
\end{pf}

\begin{corollary}
\label{pamauza}	
The abelian variety $A_0$ is contained in the abelian variety $A_1$.
\end{corollary}

\begin{pf}
Since the triple composition $r_*\circ w_*\circ \zeta _{\QQ _l}$ is $0$ by Proposition \ref{kutki}, the homomorphism $\zeta _{\QQ _l}$ factorizes through the group $H^1_{\et }(A_1,\QQ _l(1-p))$.
\end{pf}

\section{Geometric generic versus very general version of $A_0$}
\label{gengeneric}

Since the ground field $k$ is uncountable and algebraically closed, its transcendental degree over the primary subfield is infinite. We will be using the following terminology. If $S$ is an integral algebraic scheme, a (Zariski) $c$-closed subset in $S$ is a union of a countable collection of Zariski closed irreducible subsets in $S$. A (Zariski) $c$-open subset in $S$ is the complement to a $c$-closed subset in $S$. A property $\mathbf P$ of points in $S$ holds for a very general point on $S$ if there exists a $c$-open subset $U$ in $S$, such that $\mathbf P$ holds for each closed point in $U$.

The purpose of this section is to convince the reader that, given a flat family $\bcX \to S$ over an integral base $S$ over $k$, there exists a natural $c$-open subset $U$ in $S$, such that the fibres $\bcX _P$, for all closed $k$-points $P\in U$, are isomorphic to the geometric generic fibre $\bcX _{\bar \eta }$, as schemes over $\Spec (\QQ )$, and these isomorphisms preserve algebraic and rational equivalence of algebraic cycles. This is certainly a folklore, but we give all proofs for completeness. Then we use such isomorphisms to show that the abelian variety $A_0$ for fibres in a family is of purely algebraic nature, and therefore its very general and geometric generic versions coincide.

Let $S$ be an integral affine scheme of finite type over $k$, let $k(S)$ be the function field of $S$, and let $\overline {k(S)}$ be the algebraic closure of the field $k(S)$. Let $I(S)$ be the ideal of $S$ and let $f_1,\dots ,f_n$ be a set of generators of $I(S)$. The polynomials $f_i$ have a finite number of coefficients. Then we can choose a countable algebraically closed subfield $k_0$ in $k$, such that there exists an irreducible quasi-projective scheme $S_0$ over $k_0$ with $S=S_0\times _{\Spec (k_0)}\Spec (k)$.

Let $Z$ be a closed subscheme of $S_0$, and let $i_Z:Z\subset S_0$ be the corresponding closed embedding. Then $Z$ is defined by an ideal $\goa $ in $k_0[S_0]$. Since the field $k_0$ is countable and $\goa $ is finitely generated, there exists only countably many closed subschemes $Z$ in $S_0$. For each $Z$ let $U_Z$ be the complement $S_0\smallsetminus \im (i_Z)$, $Z_{k }=Z\times _{k_0}k$ and $(U_Z)_k=U_Z\times _{k_0}k$. If $(i_Z)_k:Z_k\to S$ is the pull-back of $i_Z$ with respect to the extension $k/k_0$, then $(U_Z)_k$ is the complement $S\smallsetminus \im ((i_Z)_k)$. Let, furthermore,
  $$
  U=S\smallsetminus \cup _Z\im (({i_Z})_k)=
  \cap _Z(U_Z)_{k }\; ,
  $$
where the union is taken over closed subschemes $Z$, such that $\im (({i_Z})_k)\neq S$. Notice that the last condition is equivalent to the condition $\im (i_Z)\neq S_0$. The set $U$ is $c$-open by its construction (see also the proof of Lemma 2.1 in \cite{Vial}).

\begin{lemma}
\label{iso}
For any closed $k$-point $P$ in $U$, one can construct a field isomorphism between $\overline {k(S)}$ and $k$, whose value at $f\in k_0[S_0]$ is $f(P)$.
\end{lemma}

\begin{pf}
If now $P$ is a closed $k$-point in the above defined subset $U$ in $S$, defined by the corresponding morphism $f_P:\Spec (k )\to S$, then its image under the projection $\pi :S\to S_0$ belongs to $U_Z$ for each $Z$, such that $\im (i_Z)\neq S_0$. This means that this image is noting but the generic point $\eta _0=\Spec (k_0(S_0))$ of the scheme $S_0$. In other words, there exists a morphism $h_P:\Spec (k)\to \Spec (k_0(S_0))=\eta _0$, such that $\pi \circ f_P=g_0\circ h_P$, where $g_0$ is a morphism from the generic point $\eta _0$ to $S_0$. In terms of commutative rings, it means the following. If $\ev _P:k[S]\to k$ is the evaluation at $P$, i.e. the morphism inducing $f_P$ on spectra, there exists a morphism of fields $\epsilon _P$ making the diagram
  \begin{equation}
  \label{ulitka}
  \diagram
  k[S] \ar[rr]^-{\ev _P} & & k\\ \\
  k_0[S_0] \ar[uu]_-{} \ar[rr]^-{} & & k_0(S_0)
  \ar[uu]_-{\epsilon _P}
  \enddiagram
  \end{equation}
to be commutative, where $k$ in the top right corner is considered as the residue field of the scheme $S$ at $P$. Since the left vertical morphism is injective, $k _0[S_0]\smallsetminus \{ 0\} $ is a multiplicative system in $k[S]$. It is not hard to see that the localization $(k_0[S_0]\smallsetminus \{ 0\} )^{-1}k[S]$ is the tensor product of $k[S]$ and $k_0(S_0)$ over $k_0[S_0]$. This is why there exists a unique universal morphism of rings $\varepsilon _P$ whose restriction on the ring $k[S]$ is $\ev _P$ and the restriction on $k_0(S_0)$ is $\epsilon _P$. Our aim is to construct an embedding of $k(S)$ into $k$ whose restriction to $k_0(S_0)$ would be $\epsilon _P$. Certainly, such an embedding will {\it not} be over the ring $(k_0[S_0]\smallsetminus \{ 0\} )^{-1}k[S]$.

Let $d$ be the dimension of $S_0$. By the Noether normalization lemma, there exist $d$ algebraically independent elements, $x_1,\dots ,x_d$ in $k_0[S_0]$, such that the latest ring is integral (i.e. finitely generated) over the ring $k _0[x_1,\dots ,x_d]$, and $k_0(S_0)$ is algebraic over the field of fractions $k_0(x_1,\dots ,x_d)$.
Then $k[S]$ is integral over the ring $k[x_1,\dots ,x_d]$ and $k (S)$ is algebraic over $k(x_1,\dots ,x_d)$. Let $b_i=\ev _P(x_i)$ for $i=1,\dots ,d$. Since $P\in U$, the quantities $b_1,\dots ,b_d$ are algebraically independent over $k_0$. Extend the set $\{ b_1,\dots ,b_d\} $ to a transcendental basis $B$ of $k$ over $k_0$, so that $k=k_0(B)$. Since $B$ is of infinite cardinality, so is the set $B\smallsetminus \{ b_1,\dots ,b_d\} $.

Choose and fix a bijection
  $$
  B\stackrel{\sim }{\to }
  B\smallsetminus \{ b_1,\dots ,b_d\} \; .
  $$
It gives a field embedding
  $$
  k=k_0(B)\simeq
  k_0(B\smallsetminus \{ b_1,\dots ,b_d\} )\subset
  k_0(B)=\Omega
  $$
over $k_0$, such that the set $\{ b_1,\dots ,b_d\} $ is algebraically independent over its image. The latter embedding induces a new field embedding
  $$
  k (x_1,\dots ,x_d)\to k
  $$
sending $x_i$ to $b_i$ for each $i$. The restriction of this field embedding on $k_0(x_1,\dots ,x_d)$ is the restriction of $\epsilon _P$ on the same field. Since $k (S)$ is the tensor product of $k(x_1,\dots ,x_d)$ and $k_0(S_0)$ over $k_0(x_1,\dots ,x_d)$, we get a uniquely defined embedding
  $$
  k(S)\to k\; ,
  $$
which can be extended to an isomorphism
  $$
  e_P:\overline {k(S)}\stackrel{\sim }{\lra }k\; .
  $$
As the square (\ref{ulitka}) is commutative, $e_P(f)=f(P)$ for each $f$ in $k_0[S_0]$.
\end{pf}

\begin{remark}
{\rm It is important to mention that the above isomorphism $e_P$ is non-canonical depends on the choice of the transcendental basis $B$ containing the quantities $b_1,\dots ,b_d$.}
\end{remark}

Let now $f:\bcX \to S$ be a smooth morphism of schemes over $k $. Extending $k_0$ if necessary, we may assume that there exists a morphism of schemes $f_0:\bcX _0\to S_0$ over $k_0$, such that $f$ is the pull-back of $f_0$ under the field extension from $k_0$ to $k$. Let $\eta _0=\Spec (k_0(S_0))$ be the generic point of the scheme $S_0$, $\eta =\Spec (k(S))$ the generic point of the scheme $S$, and $\bar \eta =\Spec (\overline {k(S)})$ be the geometric generic point of $S$. Then we also have the corresponding fibres $\bcX _{0,\eta _0}$, $\bcX _{\eta }$ and $\bcX _{\bar \eta }$.

Pulling back the scheme-theoretic isomorphism $\Spec (e_P)$ to the fibres of the family $f$, we obtain the Cartesian square
  $$
  \diagram
  \bcX _P \ar[dd]_-{\varkappa _P} \ar[rr]^-{} & &
  \Spec (k) \ar[dd]^-{\Spec (e_P)} \\ \\
  \bcX _{\bar \eta } \ar[rr]^-{} & & \bar \eta
  \enddiagram
  $$
Since $\Spec (e_P)$ is an isomorphism of schemes over $\eta _0$, the morphism $\varkappa _P$ is an isomorphism of schemes over $\bcX _{0,\eta _0}$.

For any field $F$, a scheme $Y$ over $F$ and an automorphism $\sigma $ of $F$ let $Y_{\sigma }$ be the fibred product of $Y$ and $\Spec (F)$ over $\Spec (F)$, with regard to the automorphism $\Spec (\sigma )$, and let $w_{\sigma }:Y_{\sigma }\stackrel{\sim}{\to }Y$ be the corresponding isomorphism of schemes over $\Spec (F^{\sigma })$, where $F^{\sigma }$ is the subfield of $\sigma $-invariant elements in $F$.

Let $L$ be a field subextension of $k/k_0$. The projection $\bcX \to \bcX _0$ naturally factorizes through ${\bcX _0}_L=\bcX _0\times _{\Spec (k_0)}\Spec (L)$. Composing the embedding of the fibre $\bcX _P$ into the total scheme $\bcX $ with the morphism $\bcX \to {\bcX _0}_L$ we can consider $\bcX _P$ as a scheme over ${\bcX _0}_L$.

If now $P'$ is another closed $k$-point in $U$, let $\sigma _{PP'}=e_{P'}\circ e_P^{-1}$ be the automorphism of the field $k$, and let $\varkappa _{PP'}=\varkappa _{P'}^{-1}\circ \varkappa _P$ be the induced isomorphism of the fibres as schemes over $\Spec (k^{\sigma _{PP'}})$. In these terms, $(\bcX _P)_{\sigma _{PP'}}=\bcX _{P'}$, the isomorphism $w_{\sigma _{PP'}}:\bcX _{P'}\stackrel{\sim }{\to }\bcX _P$ is over $\bcX _0\times _{\Spec (k_0)}\Spec (k^{\sigma _{PP'}})$, and $w_{\sigma _{PP'}}=\varkappa _{P'P}$. To see that we just need to use Lemma \ref{iso} and pull-back the scheme-theoretic isomorphisms between points on $S$ to isomorphisms between the corresponding fibres of the morphism $f:\bcX \to S$.

\begin{remark}
{\rm The assumption that $S$ is affine is not essential, of course. We can always cover $S$ by open affine subschemes, construct the system of isomorphisms $\varkappa $ in each affine chart and then construct ``transition isomorphisms" between very general fibres in a flat family over an arbitrary integral base $S$.}
\end{remark}

Let now $S$ be an integral scheme of finite type over $k$, let $\bcX $ and $\bcY $ be two schemes, both smooth, projective and connected over $S$, and let
   $$
   \xymatrix{
   \bcY \ar[rd]_-{f} \ar[rr]^-{r} & & \bcX \ar[ld]^-{g} \\
   & S & }
   $$
be a closed embedding morphism of schemes over the base $S$. Extending $k _0$ appropriately we may assume that there exist models $f_0$, $g_0$ and $r_0$ over $k_0$ of the morphisms $f$, $g$ and $r$ respectively, such that $g_0\circ r_0=f_0$. Then, for any closed $k$-point $P$ in $U$, the diagram
  \begin{equation}
  \label{ulitka3}
  \diagram
  \bcY _P \ar[dd]_-{\varkappa _P} \ar[rr]^-{r_P} & &
  \bcX _P \ar[dd]^-{\varkappa _P} \\ \\
  \bcY _{\bar \eta } \ar[rr]^-{r_{\bar \eta }}
  & & \bcX _{\bar \eta }
  \enddiagram
  \end{equation}
is commutative, where $r_P$ and $r_{\bar \eta }$ are the obvious morphisms on fibres induced by the morphism $r$. Then, of course, the isomorphisms $\varkappa _{PP'}$ commute with the morphisms $r_P$ and $r_{P'}$, for any two closed $k$-points $P$ and $P'$ in $U$. Cutting out more Zariski closed subsets from $U$ we may assume that the fibres of the families $f$ and $g$ over the points from $U$ are smooth.

\begin{lemma}
\label{ratkappa}
The scheme-theoretic isomorphisms $\varkappa _P$ preserve the algebraic and rational equivalence of algebraic cycles.
\end{lemma}

\begin{pf}
As we already explained in Section \ref{Chowmonoids}, if $\alpha :k\stackrel{\sim }{\to }k'$ is an isomorphism of fields, the functorial bijections $\Upsilon $ from the representation of Chow monoids by Chow schemes commute through the isomorphisms of monoids and Hom-sets induced by the isomorphism $\Spec (\alpha )$. In particular, if $k'=\overline {k(S)}$ and $\alpha =e_P^{-1}$, the bijections $\Upsilon (\bcX _{\bar \eta })$ over $\overline {k(S)}$ commute with the bijections $\Upsilon ({\bcX _P})$ over $k$, and the same for $\bcY $. The commutativity for the sections of the corresponding pre-sheaves on an algebraic curve $C$ over $k$ and its pull-back $C'$ over $k'$ gives the first assertion of the lemma. If $C=\PR ^1$, we get the second one.
\end{pf}

Assume now that Assumptions $1$, $2$ and $3$ are satisfied for the geometric generic fibre $\bcY _{\bar \eta }$, and for the fibre $\bcY _P$ for each closed point $P$ in $U$. Let also
 $$
 \psi _{\bar \eta }:A^p(\bcY _{\bar \eta })\stackrel{\sim }{\to }A_{\bar \eta }
 $$
and
 $$
 \psi _P:A^p(\bcY _P)\stackrel{\sim }{\to }A_P
 $$
be the corresponding regular parametrizations. By Lemma \ref{ratkappa}, the isomorphism $\varkappa _P$ induces the push-forward isomorphism of abelian groups
  $$
  {\varkappa _P}_*:A^p(\bcY _P)\to A^p(\bcY _{\bar \eta })\; .
  $$
Let
  $$
  \kappa _P:A_P\to A_{\bar \eta }
  $$
be the composition given by the commutative diagram
  \begin{equation}
  \label{som}
  \diagram
  A_P \ar[dd]_-{\kappa _P}
  \ar[rr]^-{\psi _P^{-1}} & & A^p(\bcY _P)
  \ar[dd]^-{{\varkappa _P}_*} \\ \\
  A_{\bar \eta } & & A^p(\bcY _{\bar \eta })
  \ar[ll]_-{\psi _{\bar \eta }}
  \enddiagram
  \end{equation}

Consider the obvious commutative diagram
  \begin{equation}
  \label{som}
  \diagram
  A_P \ar[dd]_-{\kappa _P}
  & & A^p(\bcY _P) \ar[ll]_-{\psi _P}
  \ar[dd]^-{{\varkappa _P}_*} & & C^p_{d,d}(\bcY _P) \ar[ll]_-{\theta _d}
   \ar[dd]^-{C^p_{d,d}(\varkappa _P)} \\ \\
  A_{\bar \eta } & & A^p(\bcY _{\bar \eta })
  \ar[ll]_-{\psi _{\bar \eta }} & & C^p_{d,d}(\bcY _{\bar \eta }) \ar[ll]_-{\theta _d}
  \enddiagram
  \end{equation}
The top and bottom horizontal compositions in this diagram are regular morphisms of schemes over $k$ and $\overline {k(S)}$ respectively, and the vertical morphism from the right hand side is a regular morphism of schemes over $\QQ $. It follows that the homomorphism $\kappa_P:A_P\to A_{\bar \eta }$ is a regular morphism of schemes over $\QQ $ too.

Now, the commutative diagram (\ref{ulitka3}) gives the commutative diagram
  \begin{equation}
  \label{nalim}
  \diagram
  A^p(\bcY _P) \ar[dd]_-{{\varkappa _P}_*}
  \ar[rr]^-{{r_P}_*} & &
  A^{p+e}(\bcX _P)\ar[dd]^-{{\varkappa _P}_*} \\ \\
  A^p(\bcY _{\bar \eta })
  \ar[rr]^-{{r_{\bar \eta }}_*}
  & & A^{p+e}(\bcX _{\bar \eta })
  \enddiagram
  \end{equation}
where $e$ is the codimension of $\bcY _{\bar \eta }$ in $\bcX _{\bar \eta }$. Let $A_{P,1}$ and $A_{\bar \eta ,1}$ be the abelian subvarieties in $A_P$ and $A_{\bar \eta }$ respectively, constructed in Section \ref{setting&assumptions}. Let, furthermore, $A_{P,0}$ and $A_{\bar \eta ,0}$ be the abelian subvarieties in $A_{P,1}$ and $A_{\bar \eta ,1}$ respectively, provided by Theorem \ref{kot} and Corollary \ref{pamauza}.

\begin{proposition}
\label{soglasovannost}
For any closed point $P$ in $U$,
  $$
  \kappa _P(A_{P,1})=A_{\bar \eta ,1}
  $$
and
 $$
 \kappa _P(A_{P,0})=A_{\bar \eta ,0}
 $$
\end{proposition}

\begin{pf}
The first claim is actually true for any closed point $P$ on $S$, not only on $U$, and can be easily deduced using specialization isomorphisms on \'etale cohomology groups. Let us prove the second claim. Let $\Xi _P$ be the countable subset in $A_P$ and $\Xi _{\bar \eta }$ the countable subset in $A_{\bar \eta }$, such that we have the kernels
  $$
  K_P=\cup _{x\in \Xi _P}(s+A_{P,0})
  \qand
  K_{\bar \eta }=\cup _{x\in \Xi _{\bar \eta }}(x+A_{\bar \eta ,0})
  $$
 in $A_P$ and $A_{\bar \eta }$ respectively (see Theorem \ref{kot}). Then
  $$
  \kappa _P(K_P)=\kappa _P(\cup _{x\in \Xi _P}(x+A_{P,0}))=
  \cup _{x\in \Xi _P}(\kappa _P(x)+\kappa _P(A_{P,0}))
  $$
The definition of $\kappa _P$ and the commutative diagram (\ref{nalim}) give that
  $$
  \kappa _P(K_P)=K_{\bar \eta }\; .
  $$
Therefore,
  $$
  \cup _{x\in \Xi _P}(\kappa _P(x)+\kappa _P(A_{P,0}))=
  \cup _{x\in \Xi _{\bar \eta }}(x+A_{\bar \eta ,0})
  $$
inside the abelian variety $A_{\bar \eta }$. Since the group isomorphisms $\kappa _P$ are regular morphisms of schemes over $\Spec (\QQ )$, we obtain that $\kappa _P(A_{P,0})$ is a Zariski closed subset in $A_{\bar \eta }$. Since $\kappa _P(A_{P,0})$ is a subgroup in $A_{\bar \eta }$, it is an abelian subvariety in $A_{\bar \eta }$. Lemma \ref{belka} and Lemma \ref{lisa} finish the proof.
\end{pf}

\begin{remark}
{\rm Of course, the set $U$ is not uniquely defined, and all the same works well over the union of all such $c$-open subset in the integral scheme $S$. The behaviour of $A_0$ beyond the union of the sets $U$ is an open question of particular importance and deserves a separate big research programme.
}
\end{remark}

\section{\'Etale monodromy argument for cycles of dimension $p-1$}
\label{hyperplane}

This is the main section of the paper, in which we apply Theorem A in the family of hyperplane sections of a projective variety embedded into a projective space. In such a case one can enhance the study of $A_0$ by the monodromy argument in terms of \'etale $l$-adic cohomology over $k$.

Let $d=2p$ and let $X$ be a smooth $d$-dimensional projective variety over the ground field $k$. Fix a closed embedding $X\subset \PR ^m$, such that $X$ is not contained in a smaller linear subspace in $\PR ^m$. Let
  $$
  \bcH =\{ (P,H)\in \PR ^m\times {\PR ^m}^{\vee }\; |\; P\in H\}
  $$
be the universal hyperplane, and let $p_1$ and $p_2$ be the projections of $\bcH $ on $\PR ^m$ and ${\PR ^m}^{\vee }$ respectively. Let
  $$
  \bcX =X\times {\PR ^m}^{\vee }
  $$
and let
  $$
  \bcY =\bcX \cap \bcH
  $$
inside $\PR ^m\times {\PR ^m}^{\vee }$.

Let
  $$
  f:\bcY \to {\PR ^m}^{\vee }
  $$
be the composition of the closed embedding $\bcY \subset \bcH $ with the projection $p_2$, let
  $$
  g:\bcX \to {\PR ^m}^{\vee }
  $$
be the composition of the closed embedding of $\bcX $ into $\PR ^m\times {\PR ^m}^{\vee }$ with the projection onto ${\PR ^m}^{\vee }$, and let
   $$
   \xymatrix{
   \bcY \ar[rd]_-{f} \ar[rr]^-{r} & & \bcX \ar[ld]^-{g} \\
   & {\PR ^m}^{\vee } & }
   $$
be the obvious closed embedding over the dual projective space.

For any scheme $S$ and for any morphism of schemes
  $$
  S\to {\PR ^m}^{\vee }
  $$
let $\bcH _S\to S$ be the pull-back of $p_2$ with respect to the morphism $S\to {\PR ^m}^{\vee }$, let $\bcY _S$ be the fibred product of $\bcY $ and $\bcH _S$ over the universal hyperplane $\bcH $, and let
  $$
  f_S:\bcY _S\to S
  $$
be the induced projection, i.e. the composition of the closed embedding of $\bcY _S$ into $\bcH _S$ and the morphism $\bcH _S\to S$. Let also $\bcX _S\to S$ be the pull-back of trivial family $\bcX \to {\PR ^m}^{\vee }$ with respect to the morphism $S\to {\PR ^m}^{\vee }$. Then we obtain the closed embedding
   $$
   \xymatrix{
   \bcY _S\ar[rd]_-{f_S} \ar[rr]^-{r_S} & & \bcX _S \ar[ld]^-{g_S} \\
   & S & }
   $$
over $S$.

Assume that the scheme $S$ is integral. Let $k(S)$ be the function field of $S$,
  $$
  \eta =\Spec (k(S))
  $$
be the generic point of $S$, $\overline {k(S)}$ be the algebraic closure of $k(S)$ and
  $$
  \bar \eta =\Spec (\overline {k(S)})
  $$
be the geometric generic point of $S$. Then we also have the closed embeddings $r_{\eta }$ and $r_{\bar \eta }$ over $\eta $ and $\bar \eta $ respectively.

As in the previous section, choose an appropriate $c$-open subset $U$ in $S$, such that the point $\bar \eta $ scheme-theoretically is isomorphic to each closed point in $P$, and assume that Assumptions $1$, $2$ and $3$ are satisfied for the geometric generic fibre
  $$
  Y_{\bar \eta }=\bcY _{\bar \eta }\; ,
  $$
and for the fibre
  $$
  Y_P=\bcY _P\; ,
  $$
for each closed point $P$ in $U$. Let also
  $$
  \psi _{\bar \eta }:A^p(Y_{\bar \eta })
  \stackrel{\sim }{\to }A_{\bar \eta }
  $$
and
  $$
  \psi _P:A^p(Y_P)\stackrel{\sim }{\to }A_P
  $$
be the corresponding regular parametrizations. Then we have the abelian subvarieties
  $$
  A_{\bar \eta,0}\subset A_{\bar \eta ,1}\subset A_{\bar \eta }
  \qqand
  A_{P,0}\subset A_{P,1}\subset A_P
  $$
for each closed point $P$ in $U$.

Let first $S=D$ be a projective line inside the dual space ${\PR ^m}^{\vee }$, such that the morphism $f_D$ is a Lefschetz pencil for the variety $X$. Let $L$ be the minimal subextension of $k(D)$ in $\overline {k(D)}$, such that the abelian varieties $A_{\bar \eta ,0}$, $A_{\bar \eta ,1}$ and $A_{\bar \eta }$ are defined over $L$. Then $L$ is finitely generated and algebraic of finite degree $n$ over $k(D)$. Let $D'$ be an algebraic curve, such that $L=k(D')$ and the embedding of $k$ into $k(D)$ is induced by a generically of degree $n$ morphism from $D'$ onto $D$. Since the closed embedding of $A_{\bar \eta ,0}$ into $A_{\bar \eta ,1}$ and  $A_{\bar \eta ,1}$ into $A_{\bar \eta }$ are now defined over $L$, there exist a Zariski open subset $U'$ in $D'$, spreads $\bcA _{\bar \eta ,0}$, $\bcA _{\bar \eta ,1}$ and $\bcA _{\bar \eta }$ of $A_{\bar \eta ,0}$, $A_{\bar \eta ,1}$ and $A_{\bar \eta }$ respectively over $U'$, and morphisms
  $$
  \bcA _{\bar \eta ,0}\to \bcA _{\bar \eta ,1}
  $$
and
  $$
  \bcA _{\bar \eta ,1}\to \bcA _{\bar \eta }
  $$
over $U'$, such that, when passing to the fibres at the geometric generic point $\bar \eta $, we obtain the closed embeddings
  $$
  A_{\bar \eta ,0}\to A_{\bar \eta ,1}
  $$
and
  $$
  A_{\bar \eta ,1}\to A_{\bar \eta }
  $$
over $\overline {k(D)}$.

Let $\alpha $ be the morphism from $\bcA $ onto $U'$, and let $\alpha _0$ and $\alpha _1$ be the morphism from $\bcA _0$ and, respectively, $\bcA _1$ onto $U'$. Since $\bcA $ is a spread of $A_{\bar \eta }$ over $U'$ and $A_{\bar \eta }$ is a projective variety over $L=k(D')$, the morphism $\alpha $ is locally projective and, therefore, proper. Similarly, the morphisms $\alpha _0$ and $\alpha _1$ are proper. Cutting more points from $D'$ we may assume that the morphisms $\alpha $, $\alpha _0$ and $\alpha _1$ are all smooth over $U'$.

Let $\eta '$ be the generic point of $D'$, let $\bar \eta '=\bar \eta $ be the geometric generic point of $D'$, let $\pi _1(U',\bar \eta )$ be the \'etale fundamental group of $D'$ pointed at $\bar \eta $, and let $\pi _1^{\tame }(U',\bar \eta )$ be the corresponding tame fundamental group. For any scheme $V$ and non-negative integer $n$ let $(\ZZ /l^n)_V$ be the constant sheaf on $V$ associated to the group $\ZZ /l^n$.

Since the morphisms $\alpha _0$, $\alpha _1$ and $\alpha $ are smooth and proper, the higher direct images
  $$
  R^1{\alpha _0}_*(\ZZ /l^n)_{\bcA _0}\; ,\; \; \;
  R^1{\alpha _1}_*(\ZZ /l^n)_{\bcA _1}
  \qqand R^1\alpha _*(\ZZ /l^n)_{\bcA }
  $$
are locally constant by Theorem 8.9, Ch. I in \cite{FreitagKiehl}. Then the fibres of these sheaves at the geometric generic point $\bar \eta $ are finite continuous $\pi _1(U',\bar \eta )$-modules, see Proposition A I.7 in loc. cit. The proper base change (see, for example, Theorem $6.1'$ on page 62 in loc. cit.) gives that
  $$
  (R^1{\alpha _0}_*(\ZZ /l^n)_{\bcA _0})_{\bar \eta }=
  H^1_{\et }({\bcA _0}_{\bar \eta },\ZZ /l^n)\; ,
  $$
  $$
  (R^1{\alpha _1}_*(\ZZ /l^n)_{\bcA _1})_{\bar \eta }=
  H^1_{\et }({\bcA _1}_{\bar \eta },\ZZ /l^n)
  $$
and
  $$
  (R^1\alpha _*(\ZZ /l^n)_{\bcA })_{\bar \eta }=
  H^1_{\et }(\bcA _{\bar \eta },\ZZ /l^n)\; .
  $$
Then we obtain that $\pi _1(U',\bar \eta )$ acts continuously on
  $$
  H^1_{\et }({\bcA _0}_{\bar \eta },\ZZ /l^n)\; ,
  $$
  $$
  H^1_{\et }({\bcA _1}_{\bar \eta },\ZZ /l^n)
  $$
and
  $$
  H^1_{\et }(\bcA _{\bar \eta },\ZZ /l^n)\; .
  $$
Passing to limits on $n$ and then tensoring with $\QQ _l$ we then obtain that $\pi _1(U',\bar \eta )$ acts continuously on

  $$
  H^1_{\et }({\bcA _0}_{\bar \eta },\QQ _l)=H^1_{\et }(A_{\bar \eta ,0},\QQ _l)\; ,
  $$
  $$
  H^1_{\et }({\bcA _1}_{\bar \eta },\QQ _l)=H^1_{\et }(A_{\bar \eta ,1},\QQ _l)
  $$
and
  $$
  H^1_{\et }(\bcA _{\bar \eta },\QQ _l)=H^1_{\et }(A_{\bar \eta },\QQ _l)\; .
  $$

The homomorphism
  $$
  \zeta _{\QQ _l}:H^1_{\et }(A_{\bar \eta ,0},\QQ _l)\to H^1_{\et }(A_{\bar \eta },\QQ _l)
  $$
is the composition of the obvious homomorphisms
  $$
  \zeta _{\QQ _l}':H^1_{\et }(A_{\bar \eta ,0},\QQ _l)\to
  H^1_{\et }(A_{\bar \eta ,1},\QQ _l)
  $$
and
  $$
  \zeta _{\QQ _l}'':H^1_{\et }(A_{\bar \eta ,1},\QQ _l)\to
  H^1_{\et }(A_{\bar \eta },\QQ _l)\; .
  $$
The action of $\pi _1(U',\bar \eta )$ naturally commutes with both $\zeta _{\QQ _l}'$ and $\zeta _{\QQ _l}''$.

Without loss of generality, we may assume that $U'$ is the pre-image of a Zariski open subset $U$ in $D$ and all the fibres of the Lefschetz pencil
  $$
  f_D:\bcY _D\to D
  $$
over the closed points of $U$ are smooth. Let
  $$
  f_{D'}:\bcY _{D'}\to D'
  $$
be the pull-back of the pencil $f_D$ with respect to the morphism $D'\to D$, and let
  $$
  f_{U'}:\bcY _{U'}\to U'
  $$
be the pull-back of $f_{D'}$ with respect to the open embedding of $U'$ to $D'$. Applying the same arguments to the morphism $f_{U'}$, we obtain the continuous action of the \'etale fundamental group $\pi _1(U',\bar \eta )$ on the group $H^{2p-1}_{\et }(Y_{\bar \eta },\QQ _l)$, and it is well known that this action is tame, in the sense that it factorizes through the surjective homomorphism from $\pi _1(U',\bar \eta )$ onto $\pi _1^{\tame }(U',\bar \eta )$.

%

For each closed point $s$ in the complement $D\smallsetminus U$ let
  $$
  \delta _s\in H^{2p-1}_{\et }(Y_{\bar \eta },\QQ _l)
  $$
be the unique up-to conjugation vanishing cycle corresponding to the point $s$ in the standard sense (see Theorem 7.1 on page 247 in \cite{FreitagKiehl}), and let 
  $$
  E\subset H^{2p-1}_{\et }(Y_{\bar \eta },\QQ _l)
  $$ 
be the $\QQ _l$-vector subspace generated by all the elements $\delta _s$, $s\in D\smallsetminus U$. In other words, $E$ is the space of vanishing cycles in $H^{2p-1}_{\et }(Y_{\bar \eta },\QQ _l)$. One can show that
  $$
  E=\ker (H^{2p-1}_{\et }(Y_{\bar \eta },\QQ _l)\to
  H^{2p+1}_{\et }(X_{\bar \eta },\QQ _l))\; ,
  $$
where $X_{\bar \eta }=X\times \bar \eta $, see Section 4.3 in \cite{WeilConjII}.

In what follows we will be using the \'etale $l$-adic Picard-Lefschetz formula for the monodromy action. For each $s\in D\smallsetminus U$ let
  $$
  \pi _{1,s}\subset \pi _1^{\tame }(U,\bar \eta )
  $$
be the so-called tame fundamental group at $s$, a subgroup uniquely determined by the point $s$ up to conjugation in $\pi _1^{\tame }(U,\bar \eta )$. In terms of \cite{FreitagKiehl}, $\pi _{1,s}$ is the image of the homomorphism
  $$
  \gamma _s:\hat \ZZ (1)\to \pi _1^{\tame }(U,\bar \eta )\;,
  $$
where $\hat \ZZ (1)$ is the limit of all groups $\mu _n$, and $\mu _n$ is the group of $n$-th roots of unity in the algebraically closed field $\overline {k(U)}$ whose exponential characteristic is $1$.

The tame fundamental group $\pi _1^{\tame }(U,\bar \eta )$ is generated by the subgroups $\pi _{1,s}$. If $u$ is an element in $\hat \ZZ (1)$, let $\bar u$ be the image of $u$ in $\ZZ _l(1)$. If now $v$ is an element in the $\QQ _l$-vector space $H^{2p-1}_{\et }(Y_{\bar \eta },\QQ _l)$ the Picard-Lefschetz formula says
   \begin{equation}
   \label{PicardLefschetz}
   \gamma _s(u)x=
   x\pm \bar u\langle x,\delta _s\rangle \delta _s\; .
   \end{equation}

\medskip

\begin{proposition}
\label{main1}
Under the assumptions above, either $A_{\bar \eta ,0}=0$ or $A_{\bar \eta ,0}=A_{\bar \eta ,1}$.
\end{proposition}

\medskip

\begin{pf}
By Proposition \ref{kutki} and the fact that the space $E$ of vanishing cycles coincides with the kernel of the Gysin homomorphism
  $$
  {r_{\bar \eta }}_*:
  H^{2p-1}_{\et }(Y_{\bar \eta },\QQ _l)\to
  H^{2p+1}_{\et }(X_{\bar \eta },\QQ _l)\; ,
  $$
we see that the image of the composition
  $$
  H^1_{\et }(A_{\bar \eta ,0},\QQ _l(1-p))
  \stackrel{\zeta _{\QQ _l}}{\lra }
  H^1_{\et }(A_{\bar \eta },\QQ _l(1-p))
  \stackrel{w_*}{\lra }
  H^{2p-1}_{\et }(Y_{\bar \eta },\QQ _l)
  $$
is contained in $E$. The homomorphism $\zeta _{\QQ _l}$ is injective and compatible with the action of $\pi _1(U',\bar \eta )$. Since $p\leq 2$, the homomorphism $w_*$ is bijective, see Remark \ref{lubanskiiles}. Then
  $$
  E\simeq H^1_{\et }(A_{\bar \eta ,1},\QQ _l(1-p))
  $$
via $\zeta _{\QQ _l}''$ and $w_*$.

Since the variety $Y_{\bar \eta }$ satisfies Assumption 1, there exists a smooth projective curve $\Gamma $ and an algebraic cycle $Z$ on $\Gamma \times Y$ over $\bar \eta $, such that the cycle class $z$ of $Z$ induces the surjective homomorphism 
  $$
  z_*:A^1(\Gamma )\to A^p(Y_{\bar \eta })\; ,
  $$
whose kernel is $G$. The homomorphism
  $$
  w_*:H^1_{\et }(A_{\bar \eta },\QQ _l(1-p))\to
  H^{2p-1}_{\et }(Y_{\bar \eta },\QQ _l)
  $$
is then induced by the composition of the embedding of the curve $\Gamma $ into its Jacobian $J_{\Gamma }$ over $\bar \eta $, the quotient map from $J_{\Gamma }$ onto the abelian variety $A=J_{\Gamma }$, also over $\bar \eta $, and the homomorphism induced by the correspondence $Z$ (see Section \ref{setting&assumptions}). Spreading out the morphisms $\Gamma \to J_{\Gamma }$ and $J_{\Gamma }\to A$, as well as the cycle $X$, over a certain Zariski open subset in $D'$, we can achieve that the homomorphism $w_*$ is compatible with the action of the fundamental group $\pi _1(U',\bar \eta )$.

This gives that the composition $w_*\circ \zeta _{\QQ _l}$ is an injection of the $\pi _1(U',\bar \eta )$-module $H^1_{\et }(A_{\bar \eta ,0},\QQ _l(1-p))$ into the $\pi _1(U',\bar \eta )$-module of vanishing cycles $E$. Let 
  $$
  E_0=\im (w_*\circ \zeta _{\QQ _l})
  $$
be the image of this injection. 

Since $U'$ is finite of degree $n$ over $U$, the group $\pi _1(U',\bar \eta )$ is\label{!!!explain why} a subgroup of finite index $n$ in the \'etale fundamental group $\pi _1(U,\bar \eta )$, and the latter group acts continuously on $E$ by the standard \'etale monodromy theory. 

We are now going to use the Picard-Lefschetz formula in order to show that $E_0$ is a $\pi _1^{\tame }(U,\bar \eta )$-equivariant subspace in $E$. Obviously, it is enough to show that for each element $\gamma _s(u)$ in $\pi _1^{\tame }(U,\bar \eta )$ and any element $x$ in $E_0$ the element $\gamma _s(u)x$ is again in the space $E_0$.

Indeed, since $\langle \delta _s,\delta _s\rangle =0$, the Picard-Lefschetz formula (\ref{PicardLefschetz}) and easy induction give that
  $$
  (\gamma _s(u))^mx=x\pm m\bar u\langle x,\delta _s\rangle \delta _s\; ,
  $$
 for a natural number $m$, whence
  $$
  \bar u\langle x,\delta _s\rangle \delta _s=
  \frac{1}{m}\left((\gamma _s(u))^mx\pm x\right)\; .
  $$
When $m$ is the index of $\pi _1(U',\bar \eta )$ in $\pi _1(U,\bar \eta )$, then $(\gamma _s(u))^m$ sits in the subgroup $\pi _1(U',\bar \eta )$, so that the right hand side of the latter formula is an element of $E_0$. Applying the Picard-Lefschetz formula again, we see that $\gamma _s(u)x$ is in $E_0$.

Thus, $E_0$ is a submodule in the $\pi _1^{\tame }(U,\bar \eta )$-module $E$. Since $E$ is known to be an absolutely irreducible (see, for example, Corollary 7.4 on page 249 in \cite{FreitagKiehl}), we see that either $E_0=0$ or $E_0=E$. In the first case $H^1_{\et }(A_{\bar \eta ,0},\QQ _l)=0$, whence $A_{\bar \eta ,0}=0$. In the second case
  $$
  \zeta _{\QQ _l}':H^1_{\et }(A_{\bar \eta ,0},\QQ _l(1-p))\to
  H^1_{\et }(A_{\bar \eta ,1},\QQ _l(1-p))
  $$
is an isomorphism, whence $A_{\bar \eta ,0}=A_{\bar \eta ,1}$.
\end{pf}

Let now $T$ be the complement to the discriminant variety of $X$ in ${\PR ^m}^{\vee }$. In other words, $T$ is the set hyperplanes in $\PR ^m$ whose intersections with $X$ are smooth. Now we want to consider the global case, when the base scheme $S$ is the scheme $T$. Let again $U$ be a $c$-open subset in $T$ constructed as in Section \ref{gengeneric}. In other words, we define $U$ by removing the images of the pull-backs of all closed embeddings into the model $T_0$ of $T$ defined over the minimal field of definition of $T$. Then $U$ is a $c$-open subset in $T$, and in the dual projective space $(\PR ^m)^{\vee }$, such that, if $\xi $ is the generic point of the projective space ${\PR ^m}^{\vee }$ and $\bar \xi $ the corresponding geometric generic point, for any closed point $P\in U$ one has the isomorphism $\varkappa _P$ between $Y_P$ and $Y_{\bar \xi }$, and for any two closed points $P$ and $P'$ in $U$ one has the scheme-theoretic isomorphism $\varkappa _{PP'}$ between $Y_P$ and $Y_{P'}$, constructed in Section \ref{gengeneric}. As above, we assume that Assumptions 1, 2 and 3 are satisfied for the fibres at $\bar \xi $ and at every closed point $P$ of the set $U$.

The next theorem is Theorem B in Introduction, and it represents the main result in the paper.

\begin{theorem}
\label{main2}
In terms above, either $A_{\bar \xi ,0}=0$, in which case $A_{P,0}=0$, or $A_{\bar \xi ,0}=A_{\bar \xi ,1}$, so that $A_{P,0}=A_{P,1}$, for any closed point $P$ in $U$.
\end{theorem}

\begin{pf}
For every closed point $P$ in ${\PR ^m}^{\vee }$ let $H_P$ be the corresponding hyperplane in $\PR ^m$. Let $\Sigma $ be a Zariski closed subset in ${\PR ^m}^{\vee }$, such that for each point $P$ in the complement to $\Sigma $ in ${\PR ^m}^{\vee }$ the hyperplane $H_P$ does not contain $X$ and the scheme-theoretic intersection of $X$, and $X\cap H_P$ is either smooth or contains at most one singular point, which is double point, see Definition 1.4 and Proposition 1.5 in \cite{FreitagKiehl}, or read through Expos\'e XVII in \cite{SGA7II}. Let $G$ be the Grassmannian of lines in ${\PR ^m}^{\vee }$. There is a Zariski open subset $W$ in $G$, such that for each line $D$ in $W$ the line $D$ does not intersect $\Sigma $ and the corresponding codimension $2$ linear subspace in $\PR ^m$ intersects $X$ transversally. In other words, any line $D$ from $W$ gives rise to a Lefschetz pencil on $X$, see the top of page 180 in \cite{FreitagKiehl} or Expos\'e XVII in \cite{SGA7II}. Let $Z$ be the complement to the above $c$-open subset $U$ in ${\PR ^m}^{\vee }$. Then $Z$ is the union of a countable collection of Zariski closed irreducible subsets in ${\PR ^m}^{\vee }$, each of which is irreducible. In particular, $Z$ is $c$-closed. It follows that the condition for a line $D\in G$ to be not a subset in $Z$ is $c$-open. By Lemma \ref{zayac}, the intersection of the corresponding $c$-open subset in $G$ with $W$ is non-empty, so that we can choose a line $D$, such that $D$ gives a Lefschetz pencil $f_D:\bcY _D\to D$ and $D\cap U\neq \emptyset $. Let $P_0$ be a point in $D\cap U$ and let $\bar \eta $ be the geometric generic point of $D$. By Proposition \ref{main1}, either $A_{\bar \eta ,0}=0$ or $A_{\bar \eta ,0}=A_{\bar \eta ,1}$. Suppose $A_{\bar \eta ,0}=0$. Proposition \ref{soglasovannost}, being applied to the pencil $f_D$, gives that $A_{P_0,0}=0$. Applying the same proposition to the family $f_T$ we obtain that $A_{\bar \xi ,0}=0$ and so for each closed point $P$ in $U$ the abelian variety $A_{P,0}$ is zero. Similarly, if $A_{\bar \eta ,0}=A_{\bar \eta ,1}$ then, by Proposition \ref{soglasovannost} applied to $f_D$ we obtain that $A_{P_0,0}=A_{P_0,1}$. Applying Proposition \ref{soglasovannost} to the family $f$ we see that $A_{\bar \xi ,0}=A_{\bar \xi ,1}$ and $A_{P,0}=A_{P,1}$ for each closed point $P$ in $U$.
\end{pf}

\section{Applications in the study of $1$-cycles on $4$-dimensional varieties}

The purpose of this section is to apply Theorem \ref{main2} (Theorem B in Introduction) in the study of rational equivalence of $1$-dimensional algebraic cycles on $4$-dimensional varieties, extending Voisin's idea on page 305 in the second volume of \cite{Voisin: The Book}. We keep all the notation and assumptions of the previous section. To enhance and, at the same time, simplify the exposition, we will assume that Assumptions 1, 2 and 3 are satisfied not only for the fibres at closed points of the $c$-open subset $U\subset T$, but rather for the fibres at all closed points of the set $T$, i.e. at all smooth sections of the variety $X$ by hyperplanes in ${\PR ^m}^{\vee }$. Assume, moreover, that $p\leq 2$ and that the group $H^{2p+1}_{\et }(X,\QQ _l)$ vanishes. The latter implies that $A_{\bar \eta ,1}=A_{\bar \eta }$ for the generic point $\eta $ of $S$, $A_{\bar \xi ,1}=A_{\bar \xi }$ and $A_{P,1}=A_P$ for each closed point $P$ in $T$.

For any closed point $P\in {\PR ^m}^{\vee }$ let $\tilde Y_P$ be the resolution of singularities of the section $Y_P$. In addition to the assumptions above, we will also require that whenever the section $Y_P$ has at worst one singular point and this point is an ordinary double point, the continuous group $A^p(\tilde Y_P)$ is weakly representable.

Since the group $A^p(Y_{\bar \xi })$ is weakly representable, we can choose a smooth projective curve $C$ over $\bar \xi $ and an appropriate algebraic cycle $Z$ on $C\times Y_{\bar \xi }$, such that the induced homomorphism $Z_*$ from $A^1(C)$ to $A^p(Y_{\bar \xi })$ is surjective. Then the homomorphism $\theta ^p_d$ from $C^p_{d,d}(Y_{\bar \xi })$ to $A^p(Y_{\bar \xi })$ is surjective for big enough $d$ (see the proof of Theorem \ref{kot}). Since the group $A^p(\tilde Y_P)$ is weakly representable, whenever $Y_P$ has at worst one singular point and this point is an ordinary double point, the homomorphism $\theta ^p_d$ from $C^p_{d,d}(\tilde Y_P)$ to $A^p(\tilde Y_P)$ is surjective as well.

It is important to stress that, as we now assume that Assumptions 1, 2 and 3 are satisfied for the fibre at every closed point $P$ of $T$, accordingly we have the abelian varieties $A_{P,0}\subset A_{P,1}\subset A_P$ for every closed point $P$ of $T$. However, it does not mean that we can extend the coherence provided by Proposition \ref{soglasovannost} from fibres at closed points of $U$ to fibres at closed points of $T$. Let $T^{\znak }$ be the set of closed points in $T$ such that $P\in T^{\znak }$ if and only if $A_{P,0}$ coincides with $A_{P,1}$.

\begin{lemma}
\label{morkov}
The set $T^{\znak }$ is constructible.
\end{lemma}

\begin{pf}
Let
  $$
  \bcV =
  \{ (Z,P)\in
  \bcC ^{p+1}_d(X)\times {\PR ^m}^{\vee }\; |\; Z
  \subset H_P\}
  $$
be the incidence subvariety, where $Z\subset H_P$ means that the codimension $p+1$ algebraic cycle $Z$ of degree $d$ on $X$ is supported on the hyperplane section $X\cap H_P$ for a closed point $P$ in ${\PR ^m}^{\vee }$. Let
  $$
  v_T:\bcV _T\to T
  $$
be the corresponding pull-back of the projection to ${\PR ^m}^{\vee }$ with respect to the inclusion of $T$ into ${\PR ^m}^{\vee }$, and let
  $$
  s_T:\bcV _T\to C^{p+1}_{d,d}(X)
  $$
be the obvious morphism from $\bcV _T$ to $C^{p+1}_{d,d}(X)$.
Let
  $$
  \bcV _T^2=\bcV _T\times _T\bcV _T
  $$
be the $2$-fold fibred product of $\bcV _T$ over $T$, and the consider the corresponding morphisms
  $$
  v_T^2:\bcV _T^2\to T
  $$
and
  $$
  s_T^2:\bcV _T^2\to C^{p+1}_{d,d}(X)\; .
  $$

By Corollary \ref{maincor} we have that $(\theta ^{p+1}_d)^{-1}(0)$ is the union of a countable collection of irreducible Zariski closed subsets in $C^{p+1}_{d,d}(X)$, say
  $$
  (\theta ^{p+1}_d)^{-1}(0)=\cup _{i\in I}Z_i\; .
  $$
Let
  $$
  W_i=(s_T^2)^{-1}(Z_i)
  $$
for each $i\in I$. For any closed point $P$ in $T$ the pre-image $(v_T^2)^{-1}(P)$ is the $2$-fold product $\bcV _P^2$ of the fibre $\bcV _P$ of the morphism $v_T$ at $P$ over $\Spec (k)$. Since the homomorphism $\theta ^p_d$ from $C^p_{d,d}(Y_P)$ to $A^p(Y_P)$ is surjective, we obtain that the condition ${r_P}_*=0$ is equivalent to the condition that the fibre $\bcV _P^2$ of the morphism $v_T^2$ at $P$ is a subset of the pre-image $\cup _{i\in I}W_i$ of $0$ under the composition $\theta ^{p+1}_d\circ s_T^2$. By Lemma \ref{zayac}, this is equivalent to saying that $\bcV _P^2$ is a subset in
  $$
  W_{i_1}\cup \dots \cup W_{i_n}\; ,
  $$
for a finite collection of indices $i_1,\dots ,i_n$ in $I$. It follows that the set $T^{\znak }$ is constructible\label{!!prove it}.
\end{pf}

We need one more easy lemma about $c$-open sets over an uncountable field.

\begin{lemma}
\label{krotto}
Let $V$ be an irreducible quasi-projective variety over $k$, and let $U$ be a nonempty $c$-open subset in $V$. Then the Zariski closure of $U$ in $V$ is $V$.
\end{lemma}

\begin{pf}
Indeed, since $U$ is $c$-open, there exists a countable union $Z=\cup _{i\in I}Z_i$ of Zariski closed irreducible subsets in $V$, such that $U=V\smallsetminus Z$. Then $\bar U$ is nothing but the complement to the interior $\interior (Z)$ of the set $Z$ in $V$. Assume that that $\interior (Z)$ is nonempty. Then there exists a nonempty subset $W$ in $\interior (Z)$, which is Zariski open in $V$. By Lemma \ref{zayac}, there exists an index $i_0\in I$, such that $W$ is contained in $Z_{i_0}$. This gives that $\interior (Z_{i_0})$ of the set $Z_{i_0}$ is nonempty. This is not possible as $Z_{i_0}$ is a closed proper subset in a Zariski topological space.
\end{pf}

Now, Bloch's definition of weak representability in \cite{BlochAnExample} (see also \cite{BlochMurre}) can be done for chow groups with coefficients in $\ZZ $ and with coefficients in $\QQ $. In the latter case we will speak about rational weak representability. Keeping the assumptions made in the beginning of this section, we can now prove Theorem C in Introduction.

\begin{theorem}
\label{main3}
If the group $A^{p+1}(X)$ is not rationally weakly representable, then the kernel of the push-forward homomorphism from $A^p(Y_P)$ to $A^{p+1}(X)$ is countable, for a very general hyperplane section $\bcY _P$.
\end{theorem}

\begin{pf}
By Theorem \ref{main2}, we have that either $A_{\bar \xi ,0}=0$ or $A_{\bar \xi ,0}=A_{\bar \xi }$. Suppose the latter. By the same Theorem \ref{main2}, $A_{P,0}=A_P$ for each closed point $P$ in the $c$-open subset $U$ in $T$. On the other hand, $U$ is a subset in $T^{\znak }$, and the set $T^{\znak }$ is constructible by Lemma \ref{morkov}. Represent $U$ as the complement to a countable union $\cup _{i\in I}D_i$ of irreducible Zariski closed subsets $D_i$ in $T$, and represent $T^{\znak }$ as a countable union $\cup _{j\in J}T^{\znak }_j$, where $T^{\znak }_j$ is Zariski open in an irreducible Zariski closed subset $Z_j$ in $T$. Let $Z$ be the union $\cup _{j\in J}Z_j$ and let $W$ be the complement to $Z$ in $T$. Then $W$ is $c$-open in $T$ and $W\cap U=\emptyset $. The intersection of $W$ and $U$ is the complement to the union of all $D_i$ and $Z_j$, $i\in I$, $j\in J$, in $T$. As $U\neq \emptyset $, it follows that $D_i\neq T$ for each index $i$. Since $W\cap U=\emptyset $, by Lemma \ref{zayac}, there must exist an index $j_0\in J$, such that $Z_{j_0}=T$. This gives that $A_{P,0}=A_P$, i.e. ${r_P}_*=0$, for each closed point $P$ in the nonempty Zariski open subset $T^{\znak }_{j_0}$ in $T$.

By Lemma \ref{krotto}, the intersection of $T^{\znak }$ with $U$ is nonempty. Let
  $$
  f_D:\bcY _D\to D
  $$
be a Lefschetz pencil for $X$, such that the set-theoretic intersection of the line $D=\PR ^1$ with the set $T^{\znak }\cap U$ is nonempty. Since the group $A^{p+1}(Y_P)$ is weakly representable for each closed point $P\in T$ and $D$ passes through $U$, it follows that the group $A^{p+1}(Y_{\bar \eta })$ is weakly representable too. Let $\Gamma _{\bar \eta }$ be a smooth projective curve and $Z_{\bar \eta }$ be an algebraic cycle of codimension $1$ on $\Gamma _{\bar \eta }\times Y_{\bar \eta }$ over $\bar \eta $ implementing the weak representability of $A^{p+1}(Y_{\bar \eta })$. Let
  $$
  D'\to D
  $$
be a finite extension of the curve $D$, such that both $\Gamma _{\bar \eta }$ and $Z_{\bar \eta }$ are defined over the function field $k(D')$. Spreading out the curve $\Gamma _{\bar \eta }$ and the cycle $Z_{\bar \eta }$ into a relative curve $\bcG \to V'$ and a relative cycle $\bcZ $ on $\bcG \times _{V'}\bcY _{V'}$ over the preimage $V'$ of a certain Zariski open subset $V$ in $D$ under the map $D'\to D$, we obtain a homomorphism
  $$
  \bcZ _*:A^1(\bcG )\to A^{p+1}(\bcY _{V'})\; .
  $$
Compactifying and resolving singularities, we obtain a surface $\bcG '$, a codimension $1$ algebraic cycle $\bcZ '$ on the variety $\bcG '\times _{D'}\bcY _{D'}$ and the homomorphism
  $$
  \bcZ '_*:A^1(\bcG ')\to A^{p+1}(\bcY _{D'})\; .
  $$

Take any element $\alpha $ in the group $A^{p+1}(\bcY _{D'})$. Let $\alpha '$ be its image in $A^{p+1}(Y_{\eta '})$, and let $\bar \alpha $ be the image of $\alpha '$ in $A^{p+1}(Y_{\bar \eta })$. Take a cycle $\bar \beta \in A^1(\Gamma _{\bar \eta })$ which goes to $\bar \alpha $ under the surjective homomorphism ${Z_{\bar \eta }}_*$ from $A^1(\Gamma _{\bar \eta })$ to $A^{p+1}(Y_{\bar \eta })$, and consider a finite extension
  $$
  D''\to D'\; ,
  $$
such that $\bar \beta $ is defined over the function field $k(D'')$.

Let $\eta '$ and $\eta ''$ be the generic points of $D'$ and $D''$ respectively and consider the following commutative diagram
    $$
    \xymatrix{
    A^1(\Gamma _{\bar \eta })
    \ar[rr]^-{{Z_{\bar \eta }}_*} & &
    A^{p+1}(Y_{\bar \eta }) \\ \\
    A^1(\Gamma _{\eta ''})
    \ar[rr]^-{{Z_{\eta ''}}_*}
    \ar[uu]^-{} & &
    A^{p+1}(Y_{\eta ''})
    \ar[uu]_-{} \\ \\
    A^1(\Gamma _{\eta '})
    \ar[rr]^-{{Z_{\eta '}}_*}
    \ar[uu]^-{} & &
    A^{p+1}(Y_{\eta '})
    \ar[uu]_-{} \\ \\
    A^1(\bcG ') \ar[rr]^-{\bcZ '_*}
    \ar[uu]^-{} & &
    A^{p+1}(\bcY _{D'})
    \ar[uu]_-{}
    }
    $$
which illustrates what is going on.

The cycle class $\bar \beta $ comes from a cycle class $\beta ''\in A^1(\Gamma _{\eta ''})$ under the pullback from $A^1(\Gamma _{\eta ''})$ to $A^1(\Gamma _{\bar \eta })$. Let $\beta '$ be the image of $\beta ''$ under the pushforward homomorphism from $A^1(\Gamma _{\eta ''})$ to $A^1(\Gamma _{\eta '})$, and let $\beta $ be a cycle class in $A^1(\bcG ')$ going to $\beta '$ under the surjective homomorphism from $A^1(\bcG ')$ to $A^1(\Gamma _{\eta '})$. Let $\gamma $ (respectively, $\gamma '$ and $\gamma ''$) be the image of the cycle class $\beta $ (respectively, $\beta '$ and $\beta ''$) under the homomorphism $\bcZ '_*$ (respectively, ${Z_{\eta '}}_*$ and ${Z_{\eta ''}}_*$). Without loss of generality, we may assume that $\eta ''/\eta '$ is Galois. Let $N$ be the corresponding norm on the group $A^{p+1}(Y_{\eta ''})$. Since the kernel of the homomorphism from $A^{p+1}(Y_{\eta ''})$ to $A^{p+1}(Y_{\bar \eta })$ is torsion, there exists a positive integer $m$, such that
  $$
  m(\gamma ''-\alpha '')=0\; .
  $$
Then
  $$
  mN(\gamma '')=mN(\alpha '')=mn\alpha '\; ,
  $$
where
  $$
  n=[\eta '':\eta ']\; .
  $$
It follows that
  $$
  m\gamma -mn\alpha
  $$
belongs to the kernel of the homomorphism from $A^{p+1}(\bcY _{D'})$ to $A^{p+1}(Y_{\bar \eta })$.

We see that, if
  $$
  B_1=\im (A^1(\bcY ')\stackrel{\bcZ '_*}{\lra }
  A^{p+1}(\bcY _{D'})
  $$
and
  $$
  B_2=
  \ker (A^{p+1}(\bcY _{D'})\to A^{p+1}(Y_{\bar \eta })\; ,
  $$
the $\QQ $-vector space $A^{p+1}(\bcY _{D'})\otimes \QQ $ is generated by the vector subspaces $B_1\otimes \QQ $ and $B_2\otimes \QQ $.

Now, as $A^{p+1}(Y_{\eta '})$ is the colimit of the groups $A^{p+1}(\bcY _{W'})$, where $W'$ runs through Zariski open subsets in $D'$, the localization sequences for open embeddings $\bcY _{W'}\subset \bcY _{D'}$ show that kernel $B_2$ is generated by the image of the homomorphism
  $$
  \oplus _{P'\in D'}A^p(Y_{P'})\to
  A^{p+1}(\bcY _{D'})\; ,
  $$
induced by the proper push-forward homomorphisms ${r_{P'}}_*$. If $D^{\znak }$ is the intersection of $T^{\znak }$ and $D$, and ${D'}^{\znak }$ is the pre-image of $D^{\znak }$ under the finite map from $D'$ onto $D$, then ${r_{P'}}_*=0$ for each closed point $P'$ in ${D'}^{\znak }$. It follows that $B_2$ is generated by the image of the homomorphism
  $$
  \oplus _{P'\in D'\smallsetminus {D'}^{\znak }}
  A^p(Y_{P'})\to A^{p+1}(\bcY _{D'})\; .
  $$
Notice that the complement $D'\smallsetminus {D'}^{\znak }$ is finite.

Next, if $Y_{P'}$ is smooth, then $Y_{P'}=Y_P$, where $P$ is the image of $P'$ under the finite map from $D'$ onto $D$, and $A^p(Y_{P'})$ is isomorphic to the abelian variety $A_P$ via the universal regular morphism $\psi _P$. In particular, $A^p(Y_{P'})$ is weakly representable. If the section $Y_{P'}=Y_P$ is singular, resolving the double point on it we obtain a nonsingualr variety $\tilde Y_P$ whose group $A^p(\tilde Y_P)$ is weakly representable by our assumption.

Thus, $B_2$ is covered by the finite direct product of weakly representable groups $A^p(Y_P)$ and $A^p(\tilde Y_P)$. It means that $B_2$ itself is weakly representable. Then, of course, $B_2\otimes \QQ $ is rationally weakly representable. The image $B_1$ of the homomorphism $\bcZ '_*\otimes \QQ $ is weakly representable because $A^1(\bcY ')$ is parametrized by the Picard variety of the surface $\bcY '$. Since $A^{p+1}(\bcY _{D'})\otimes \QQ $ is generated by rationally weakly representable $\QQ $-vector subspaces $B_1$ and $B_2$, the whole group $A^{p+1}(\bcY _{D'})$ is rationally weakly representable. This contradicts to the third assumption of the theorem.

Hence, $A_{\bar \xi ,0}=0$, and Theorem \ref{main2} finishes the proof.
\end{pf}

\section{Application to hyperplane sections of cubic fourfolds in $\PR ^5$}

Let first $X$ be a $K3$-surface embedded appropriately into $\PR ^m$. Since smooth hyperplane sections of a projective surface are smooth projective curves, Assumptions 1, 2 and 3 are satisfied. By Beauville's result, \cite{Beauville}, the group $A^2(X)$ is divisible. The third cohomology of a $K3$-surface vanishes, so that the Albanese variety is trivial. By Roitman's theorem, \cite{Roitman_tors}, the group $A^2(X)$ is uniquely divisible. Then there is no difference between rational and integral weak representability for this group. Moreover, we know that $A^2(X)$ is not representable by Mumford's result in \cite{Mumford}. By Theorem \ref{main3}, for a very general hyperplane section $Y_P$ of the surface $X$ the kernel of the push-forward homomorphism ${r_P}_*$ from $A^p(Y_P)$ to $A^{p+1}(X)$ is countable. This is, of course, a particular case of Proposition 2.4 in \cite{Voisin: Sympl. involutions}.

Consider now the case when $X$ is a smooth cubic hypersurface in $\PR ^5$. For such $X$ we have that $H^5_{\et }(\bcX ,\QQ _l)=0$ and so $A_1=A$ for $Y_{\bar \xi }$ and every smooth section $Y_P$ of the fourfold $X$. Any such a section smooth is a cubic $3$-fold in $H_P\simeq \PR ^4$, whose group $A^2(Y_P)$ is well known to be representable by the corresponding Prymian variety $\Prym (Y_P)$, see \cite{Beauville}. Since the Prym construction is of purely algebraic-geometric nature, we can do it over $\bar \xi $ getting the Prym variety $\Prym (Y_{\bar \xi })$ for the geometric generic fibre $Y_{\bar \xi }$. In other words, all the Assumptions 1, 2 and 3 are satisfied for $Y_{\bar \xi }$, as well as for all smooth hyperplane sections $Y_P$.

If a hyperplane section $Y_P$ of the cubic fourfold $X$ has one singular point and this point is an ordinary double point, then the singular cubic $Y_P$ is rational, so that $\tilde Y_P$ is rational. It follows that the group $A^2(\tilde Y_P)$ is weakly representable. If $Y_P$ is smooth, then it is unirational and so rationally connected. Hence, $A^3(Y_P)$ is trivial. The group $A^3(X)$ is not weakly representable by Theorem 0.5 in \cite{Schoen}. Since, moreover, it is uniquely divisible, see Theorem 4.7(iii) in \cite{Shen_involutions}, it is also not rationally weakly representable. Thus, all the assumptions of Theorem \ref{main3} are also satisfied.

By Theorem \ref{main3}, for each closed point $P$ in the $c$-open subset $U$ of ${\PR ^m}^{\vee }$ there exists a countable set $\Xi _P$ of closed points in the Prymian $\Prym (Y_P)$ of the hyperplane section $Y_P$, such that the kernel of the homomorphism ${r_P}_*$ from $\Prym (Y_P)$ to $A^3(X)$ is countable.

In particular, if $\Sigma $ and $\Sigma '$ are two linear combinations of lines of the same degree on $X$, supported on $Y_P$, then $\Sigma $ is rationally equivalent to $\Sigma '$ on $X$ if and only if the point on $\Prym (Y_P)$, represented by the class of $\Sigma -\Sigma '$, occurs in $\Xi _P$.

Notice also that the group $A^3(Y_{\eta })$ can be nonzero, but we know that it is torsion. Since $A^2(Y_P)$ is divisible, any cycle class in $A^3(X)$ is represented, up to torsion, by line configurations supported on hyperplane sections.

\bigskip

\begin{small}

\end{small}

\bigskip

\bigskip

\begin{small}

{\sc School of Mathematics, Tata Institute of Fundamental Research, Homi Bhabha Road, Colaba, Mumbai 400005, India}

\medskip

\footnotesize{{\it E-mail addres}: {\tt kalyan484@gmail.com}}

\bigskip

{\sc Department of Mathematical Sciences, University of Liverpool, Peach Street, Liverpool L69 7ZL, England, UK}

\medskip

\footnotesize{{\it E-mail address}: {\tt vladimir.guletskii@liverpool.ac.uk}}

\end{small}

\end{document}